# Symmetry breaking perturbations and strange attractors.

Anna Litvak Hinenzon and Vered Rom-Kedar
The faculty of mathematical sciences
The Weizmann Institute of Science, P.O.B. 26, Rehovot 76100, Israel.
litvak, vered@wisdom.weizmann.ac.il

**Abstract**

The asymmetrically forced, damped Duffing oscillator is introduced as a prototype model for analyzing the homoclinic tangle of symmetric dissipative systems with *symmetry breaking* disturbances. Even a slight fixed asymmetry in the perturbation may cause a substantial change in the asymptotic behavior of the system, e.g. transitions from two sided to one sided strange attractors as the other parameters are varied. Moreover, slight asymmetries may cause substantial asymmetries in the relative size of the basins of attraction of the unforced nearly symmetric attracting regions. These changes seems to be associated with homoclinic bifurcations. Numerical evidence indicates that *strange attractors* appear near curves corresponding to specific secondary homoclinic bifurcations. These curves are found using analytical perturbational tools.

## 1 Introduction

The forced and damped Duffing oscillator:

$$\ddot{x} + \varepsilon\delta\,\dot{x} - x + x^3 = \varepsilon\gamma\cos(\omega t),\ \ (x,t) \in \Re^1 \times \Re^1, \tag{1}$$

has served as a prototype model for investigating low dimensional chaotic behavior in enumerate publications (see [16, 22, 37] and references therein). Its significance lies in its simple "typical form" which appears in many applications. Indeed, the unperturbed Duffing oscillator represents the normal form for Hamiltonian systems with $Z_2$ symmetry [13]. Thus, whether different types of perturbations lead to substantially different dynamics is of mathematical and physical significance. The perturbation of (1) has two specific properties - it has no non-linear terms in $x, \dot{x}$ and it is symmetric; (1) is invariant under $x \to -x, t \to t + \pi/\omega$.

Numerical simulations suggest that the inclusion of nonlinear dissipation term in the perturbation does not alter the qualitative behavior of the forced system [29]. Namely, no new bifurcation sequences or new types of attractors appear, though the location of the various bifurcation curves of (1) changes. The effect of asymmetric potentials has been investigated when the forcing is adiabatic, see [25, 3] and references therein. In this paper, we examine the effect of asymmetric forcing on the Duffing oscillator by introducing the *asymmetrically forced, damped, Duffing oscillator* (AFDO):

$$\ddot{x} + \delta\,\dot{x} - x + x^3 = (x - \beta x^2)\gamma\cos(\omega t),\ \ (x,t) \in \Re^1 \times \Re^1. \tag{2}$$

which contains the asymmetry perturbation parameter $\beta$. Here we show that the inclusion of the *physically typical* asymmetric forcing perturbations alters the *qualitative* behavior of the system in some range of parameter values.

Theoretical and numerical investigations of forced and damped systems with homoclinic tangle is problematic since these may attain Strange Attractors (SA)[1] and periodic sinks simultaneously. Moreover, the existence of the SA is extremely sensitive to changes in parameter values. The existence of Newhouse sinks near homoclinic tangencies implies that small changes in the parameters may destroy the SA. These observations are reflected in the difficulties of proving the existence of SA in such systems (see reviews [16, 37]). Analytical results regarding the existence of SA [1], their basins of attraction [5], and the construction of unique natural invariant measure (the SRB measure) [2], has been recently published for the Hénon map. These proofs are in the strong dissipation limit, for which the strange attractor appears as a one dimensional attractor multiplied by a cantor set. Some of these results may be applied locally to neighborhoods of homoclinic tangencies [23].

[1] These are attractors with sensitive dependence on initial conditions, i.e. attractors which have a dense orbit with positive Lyapunov exponent (see [9]).

It follows that SA are expected to emerge near homoclinic tangencies. We use analytical tools for locating primary[2] homoclinic tangencies (the Melnikov analysis, [16]) *and* secondary [3] homoclinic tangencies (the SMF [33]). See [34, 17, 18, 4] for other works on the subject of multi-pulse homoclinic orbits in other settings. This presents the first application of the SMF to a dissipative system. Since homoclinic bifurcations are considered an important source of structural instabilities of dynamical systems [27, 11, 12, 14, 19, 26], their location in parameter space should indicate regions in which dramatic structural changes appear. Clearly higher order tangencies exist as well, and finding them all is a useless mission, in particular in view of Newhouse work. The philosophy here is that not all homoclinic tangencies have the same significance: primary tangencies are more important then secondary, secondary more then third order etc.. Thus there is a sense in locating the bifurcation curves of the lower order homoclinic tangencies. This approach is backed up by the TAM (topological approximation method [31, 32]), which asserts that many features of the dynamical system are determined already by the characteristics of the primary and secondary homoclinic orbits. The TAM was developed for non-dissipative systems and has been recently generalized to dissipative systems [20].

The last part of this work consists of a numerical search for SA at parameter values which are close to the analytically predicted bifurcation curves. SA have been observed in various systems exhibiting homoclinic chaos, including the forced and damped Duffing oscillator [16], the Hénon map [19], and the forced and damped cubic potential [21]. In the latter work the correspondence between the appearance of homoclinic tangencies of specific character and SA has been noted, correspondence which seems to persist for the AFDO.

This paper is ordered as follows: In section 2 we present the basic phase space structure of the AFDO, the Melnikov analysis and the bifurcation curves for primary and secondary homoclinic bifurcations. Numerical evidence suggesting the existence of SA near specific homoclinic bifurcation curves is presented in section 3, as are the typical size and shape of the basin of attractions of the attractors. Conclusions and a discussion are presented in section 4.

[2]These are one-loop homoclinic orbits which are $O(\varepsilon)$ close to the unperturbed homoclinic orbits for $t \in (-\infty, \infty)$.

[3]These are two-loops homoclinic orbits which are $O(\varepsilon)$ close to the unperturbed homoclinic orbits for $t \in (-\infty, t_0], [t_1, \infty)$, see also section 2.

# 2 Templates of the homoclinic tangle

## 2.1 Basic properties of the AFDO

Introducing the phase space coordinates $(x, y) \in \Re^2$, one rewrites equation (2), as:

$$\begin{aligned} \dot{x} &= y, \\ \dot{y} &= x - x^3 + (x - \beta x^2)\varepsilon\gamma\cos(\omega t) - \varepsilon\delta y \ , \end{aligned} \tag{3}$$

Physically, $\delta$ represents the dissipation (the damping), $\gamma$ the amplitude of the forcing, $\omega$ the frequency and $\beta$ the asymmetry disturbances. These parameters are real and by symmetry may be taken to be non negative. $\varepsilon$ is a "perturbation scaling parameter", assumed to be small. For $\varepsilon \neq 0$, there are two differences between the AFDO (2) and the Duffing oscillator ( 1). The substantial difference is that (2) includes the asymmetry parameter, $\beta$. The second difference is that for convenience, with no loss of generality, the symmetry $x \to -x, t \to t + \pi/\omega$ of (1) is replaced by the symmetry $x \to -x$ for $\beta = 0$ in (2), hence the origin is fixed for all $\varepsilon, \beta$. $\beta \neq 0$ corresponds to symmetry breaking disturbances.

The unperturbed system corresponds to the integrable Hamiltonian system with a symmetric quartic potential:

$$V(x) = -\frac{x^2}{2} + \frac{x^4}{4} \ , \tag{4}$$

and with the Hamiltonian function (energy):

$$H(x, y) = \frac{y^2}{2} + V(x) = \frac{y^2}{2} - \frac{x^2}{2} + \frac{x^4}{4} \ . \tag{5}$$

The unperturbed system, which is identical to that of the unperturbed Duffing oscillator, has three equilibrium points: two centers at $(x, y) = (\pm 1, 0)$, and a saddle at $(x, y) = (0, 0)$. The saddle point is connected to itself by two homoclinic orbits, with periodic orbits nested within and around them. The period of the unperturbed periodic orbits, $P(H)$, has the following asymptotic expansion near $H = 0$ (exact formulae for all $H$ are available [16, 37]):

$$P(H) = \begin{cases} \ln(\frac{16}{-H})(1 + O(H)), & H \to 0- \\ \\ 2\ln(\frac{16}{H})(1 + O(H)), & H \to 0+ \ . \end{cases} \tag{6}$$

In the unperturbed system the stable and the unstable manifolds of the saddle point $(0,0)$ coincide. For $\delta > 0$, and $\gamma = 0$, the unstable manifold of the saddle point near the origin falls into the two sinks created near $(\pm 1, 0)$. As for the Duffing oscillator, it may be proved that for sufficiently small values of $\gamma$ the closure of the unstable manifold (which contains the saddle and the sinks) is an attracting set of (3).

A Poincaré map in time is used to simplify the phase space portrait for the time dependent system ($\gamma \neq 0$). Keeping $\delta > 0$, and increasing $\gamma$, the following scenario occurs on both sides of the fixed point; for small values of $\gamma$, the Poincaré map is topologically equivalent to the Poincaré map with $\gamma = 0$, which is structurally stable. As $\gamma$ increases, resonance bands of higher period and higher amplitudes are created. As $\gamma$ is further increased, in addition to the resonances, a homoclinic bifurcation occurs, after which the stable and the unstable manifolds of the saddle point of the Poincaré map intersect in *transversal homoclinic orbits.* The presence of these orbits implies the existence of a complicated non-wandering Cantor set which possesses infinitely many unstable periodic orbits of arbitrary long period as well as bounded non-periodic motions. The Smale - Birkhoff Homoclinic Theorem implies that in this case the system has chaotic dynamics. When $\beta = 0$, as in the forced Duffing oscillator (1), the sequence of bifurcations described above occurs simultaneously on both sides of the fixed point. When $\beta > 0$ this changes as described below.

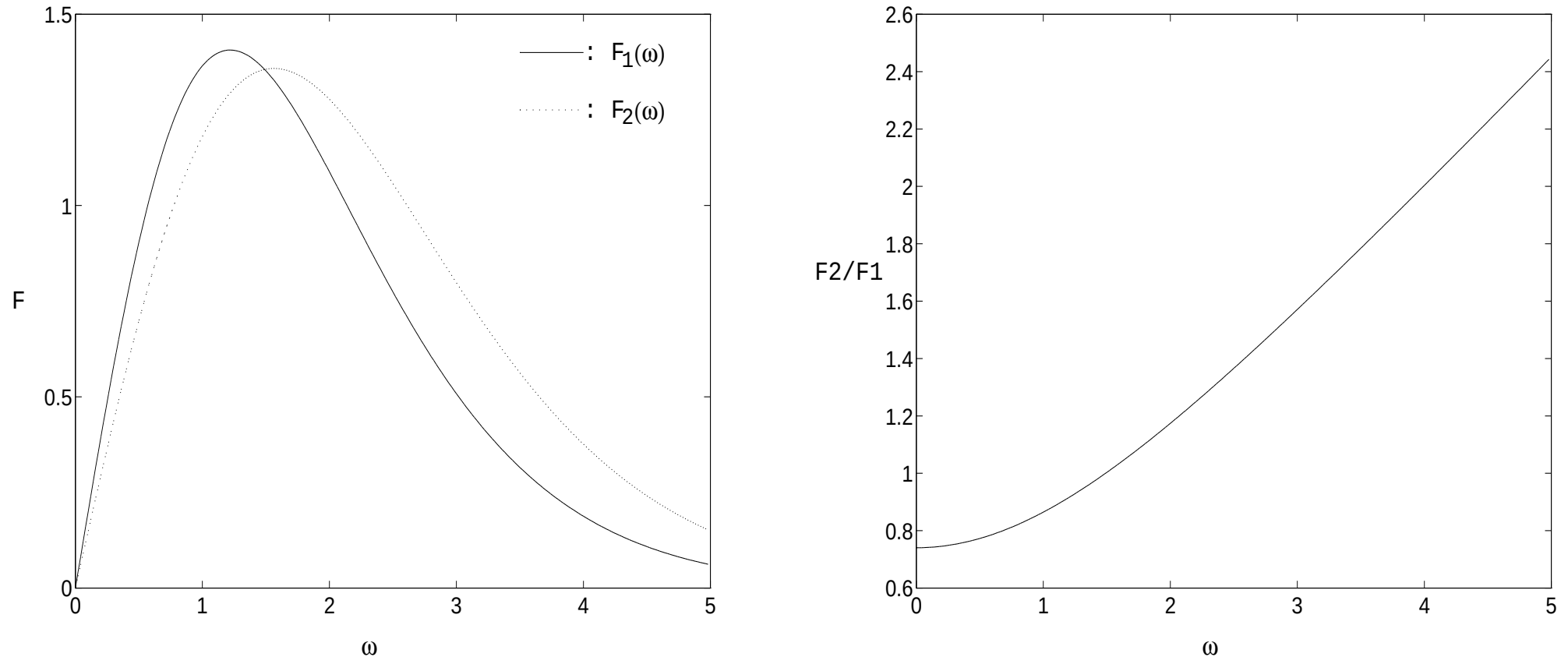

**a** **b**

a) $F_1(\omega)$, $F_2(\omega)$, b) The ratio $\frac{F_2(\omega)}{F_1(\omega)}$.

Figure 1: The Melnikov function amplitude.

The Melnikov function, $M(t_0 + \theta/\omega)$, measures the signed distance between the stable and the unstable manifolds of a hyperbolic fixed point (up to a multiplication by a constant). This distance is measured at the Poincaré section $\omega t = \theta$, and $t_0$ represents a parameterization along the unstable manifold. For the AFDO (equation (3)) two Melnikov functions are defined: $M_r(t_0; \gamma, \omega, \beta, \delta) \equiv M_r(t_0)$ (respectively $M_l(t_0)$) measures the signed distance between the right (respectively left) branches of the stable and unstable manifolds. These functions are given by:

$$\begin{aligned} M_{r,l}(t_0; \gamma, \omega, \beta, \delta) &= \int_{-\infty}^{\infty} (yx(1-\beta x)\gamma \cos(\omega t) - \delta y^2) \,|_{(q^0(t)_{r,l}, t+t_0)}\, dt \\ &= \gamma \sin(\omega t_0) F_{r,l}(\omega, \beta) - \frac{4}{3}\delta, \end{aligned} \tag{7}$$

where $q^0(t)_{r,l}$ are the right and left unperturbed homoclinic orbits of the system:

$$q^0(t)_{r,l} = \pm(\sqrt{2}\,\mathrm{sech}\, t, -\sqrt{2}\,\mathrm{sech}\, t \tanh t), \tag{8}$$

and

$$F_r(\omega, \beta) = [F_1(\omega) - \beta F_2(\omega)], \tag{9}$$

$$F_l(\omega, \beta) = [F_1(\omega) + \beta F_2(\omega)], \tag{10}$$

$$F_1(\omega) = \pi \omega^2 \,\mathrm{csch}(\frac{\pi\omega}{2}), \tag{11}$$

$$F_2(\omega) = \frac{\sqrt{2}}{3}\pi\omega(1+\omega^2)\,\mathrm{sech}(\frac{\pi\omega}{2}). \tag{12}$$

Figure 1 shows $F_1(\omega)$, $F_2(\omega)$, and the relation $\frac{F_2(\omega)}{F_1(\omega)}$. Notice that $F_1(\omega)$ and $F_2(\omega)$ are non-negative for all $\omega$.

## 2.2 Primary homoclinic intersection points

For $\varepsilon$ sufficiently small simple (respectively degenerate) zeros of the Melnikov function imply primary homoclinic transverse intersections (respectively tangencies) of the stable and the unstable manifolds of the hyperbolic fixed point [16]. Requiring $M(t_0) = M'(t_0) = 0$, it follows from (7) that primary homoclinic bifurcations occur near

$$\varepsilon\gamma_{r,l}(\delta; \omega, \beta) = \varepsilon\delta \left| \frac{4}{3[F_1(\omega) \mp \beta F_2(\omega)]} \right| \equiv \varepsilon\delta R_0^{\pm}(\omega, \beta). \tag{13}$$

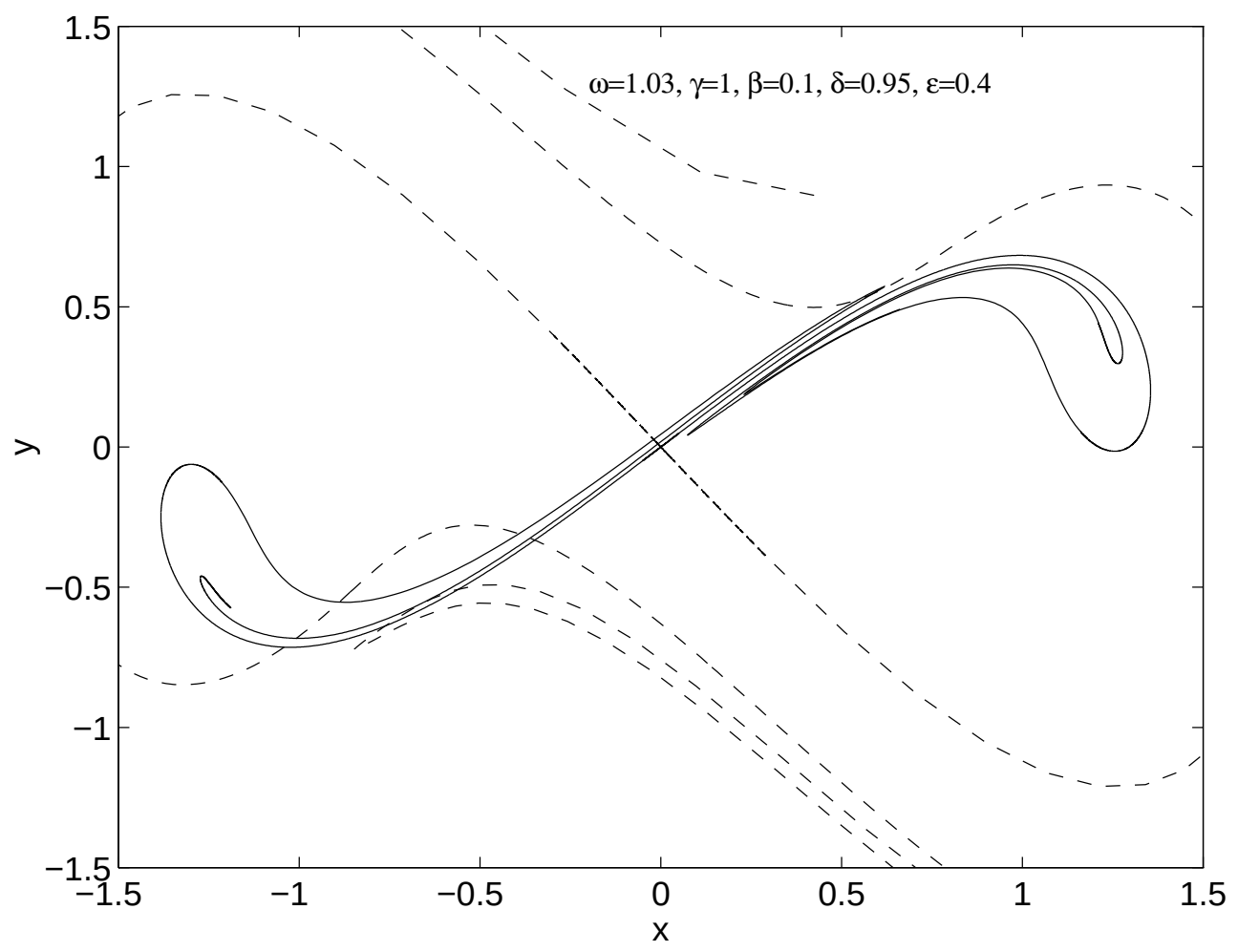

- -: stable manifold, —: unstable manifold.

Figure 2: Intersections of stable and unstable manifolds.

Hence, if

$$\frac{3\gamma\,|F_1(\omega)-\beta F_2(\omega)|}{4} < \delta \leq \frac{3\gamma[F_1(\omega)+\beta F_2(\omega)]}{4}, \tag{14}$$

then the left branches of the stable and the unstable manifolds intersect, while the right branches do not, see for example figure 2 (the intersections on the right hand side of this figure correspond to secondary homoclinic points, as described in 2.3). Similarly, if

$$\delta \leq \frac{3\gamma\,|F_1(\omega)-\beta F_2(\omega)|}{4}, \tag{15}$$

then the stable and unstable manifolds intersect on both left and right sides of the hyperbolic fixed point, see for example figure 7. It follows that the parameter space is divided to three regions (see figures 3, 4):

I For $\frac{\gamma}{\delta} < R_0^-(\omega,\beta)$ there are no primary intersections of the stable and the unstable manifolds.

II For $R_0^-(\omega,\beta) < \frac{\gamma}{\delta} < R_0^+(\omega,\beta)$ the stable and the unstable manifolds have primary intersection points on the left side of the saddle point, and do not have primary intersection points on the right side (asymmetric behavior).

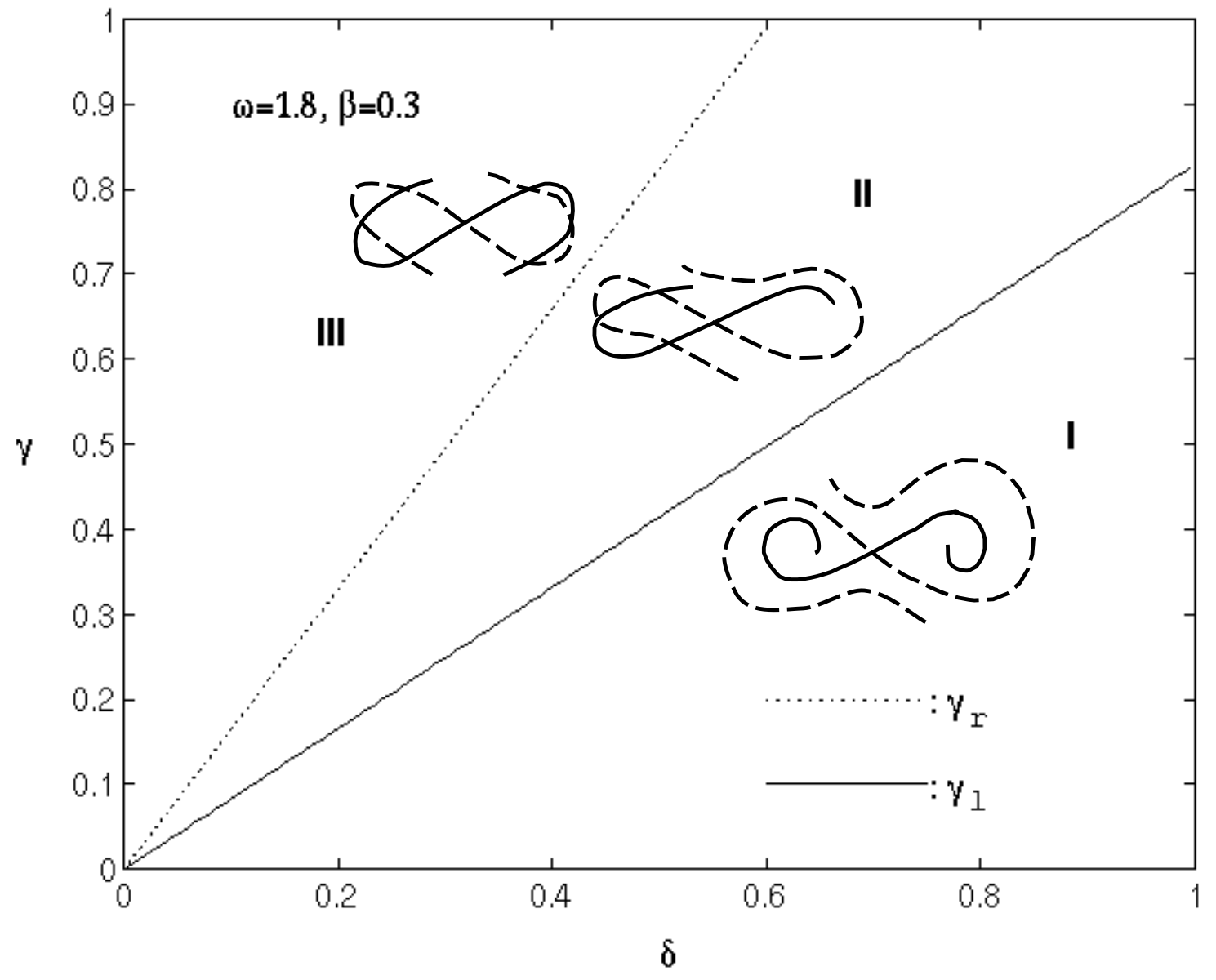

Figure 3: Primary homoclinic bifurcation curves.
Schematic phase space portraits are shown in each region.

III For $\frac{\gamma}{\delta} > R_0^+(\omega, \beta)$ primary intersections of the stable and the unstable manifolds occur both on the left and the right sides of the saddle point (like in (1), but in an asymmetric manner for $\beta \neq 0$).

Since

$$\frac{R_0^-(\omega,\beta)}{R_0^+(\omega,\beta)} = \left| \frac{1 - \beta\frac{F_2(\omega)}{F_1(\omega)}}{1 + \beta\frac{F_2(\omega)}{F_1(\omega)}} \right| \equiv r\left(\beta\frac{F_2(\omega)}{F_1(\omega)}\right), \tag{16}$$

the relative size of region II depends on the values of $x = \beta\frac{F_2(\omega)}{F_1(\omega)}$, and may be derived from the graph of $r(x)$ (figure 5) for the corresponding $\beta$ and $\omega$ values; Fixing $\beta \neq 0$, the relative size of region II varies with $\omega$ as described below. First, notice that $\frac{F_2(\omega)}{F_1(\omega)}$ is bounded from below by a positive constant $c_0 \approx 0.7$, and that it grows monotonically (in fact asymptotically linearly) with $\omega$ (see figure 1). Therefore, and since $r(x) = 1$ at $x = 0, \infty$ only, it follows that for any finite non-vanishing $\omega$ value region II is of non-vanishing measure (see for example figures 5 and 4). Since $x$ grows monotonically with $\omega$, and $x > \beta c_0$, the relative size of region II increases with $\omega$ up to the threshold value $\omega^*(\beta)$ at which $\frac{F_2(\omega^*)}{F_1(\omega^*)} = \frac{1}{\beta}$ (i.e. $x(\beta, \omega^*(\beta)) = 1)^4$. For $x > 1$, $r(x)$ strictly increases with $x$, thus for $\omega > \omega^*(\beta)$ the relative size of region II decreases as $\omega$ increases. At $x = 1$, $r(x) = 0$, namely $R_0^+(\omega, \beta) \to \infty$ as $\omega \to \omega^*$. For these values, regions I and II occupy most of the parameter space, and region III shrinks till it disappears, to order $\varepsilon$, at $\omega = \omega^*(\beta)$.

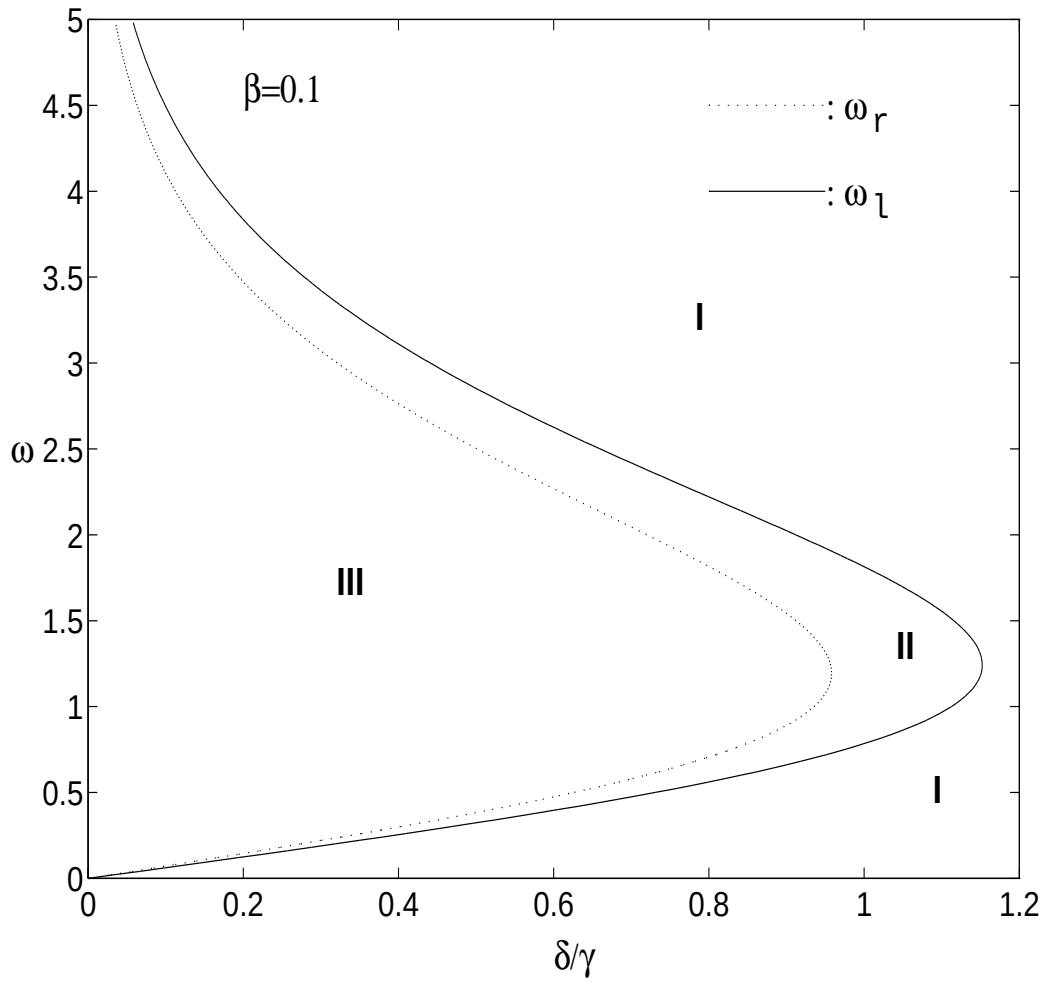

Figure 4: Primary homoclinic bifurcation curves (in $(\frac{\delta}{\gamma}, \omega)$ space).

**a** 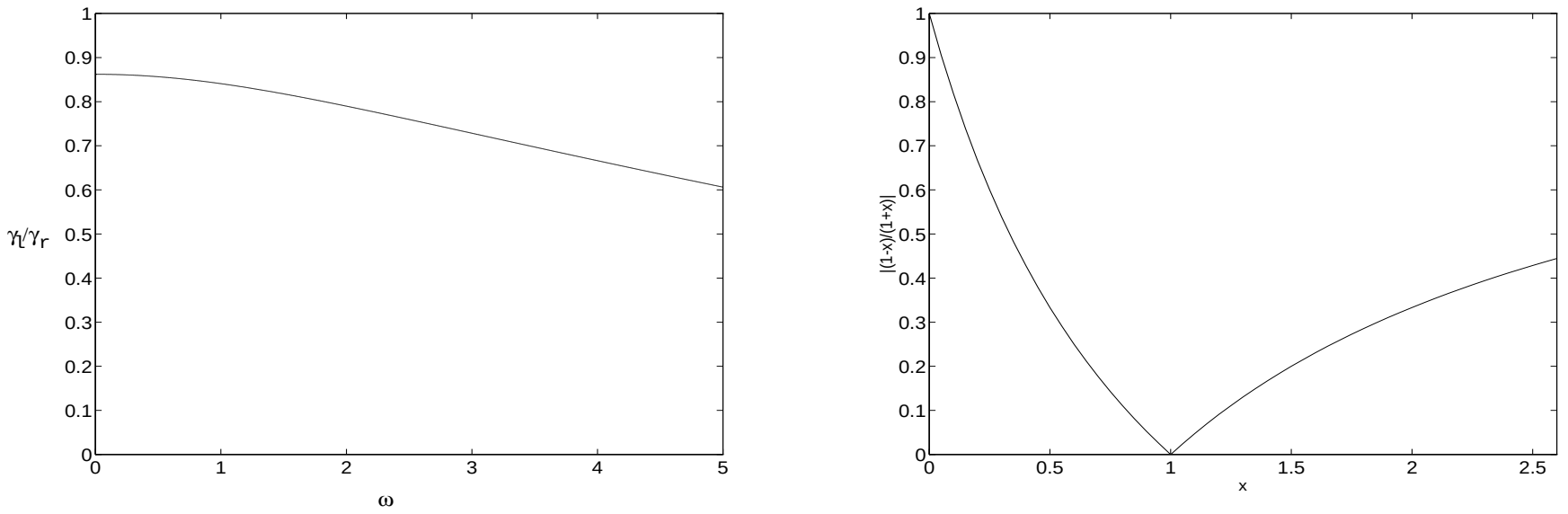 **b**

a) $\frac{\gamma_l}{\gamma_r}(\omega, \beta = 0.1)$, see eq. (13 ), b) $r(x)$, see eq. (16).

Figure 5: The relative size of region II.

[4] If $\beta > \max_\omega(\frac{F_1(\omega)}{F_2(\omega)}) = \frac{1}{c_0}$, then $\omega^*(\beta)$ does not exist, and may be considered as infinite.

## 2.3 Secondary homoclinic intersection points.

First, we describe geometrically what are the secondary homoclinic intersection points, what are their transition numbers and what are the structural indices of a homoclinic tangle. Then, we present the analytical (perturbational) method for finding these structural indices.

Consider values of the dissipation parameter $\delta$ for which the Melnikov function $M_l(t_0;\mu,\delta)$ (respectively $M_r(t_0;\mu,\delta)$) has two simple zeros (see section 2.2). Denote the corresponding PIPs (Primary homoclinic Intersection Points), ordered by the direction of the unstable manifold, by $pl_0$ (respectively $pr_0$) and $ql_0$ (respectively $qr_0$) see figure 6. Also, denote their ordered images under the Poincaré map $F$ by $pl_i, ql_i$ (respectively $pr_i, qr_i$), $i = 0, \pm1, \pm2, \pm3, ...$, i.e. $F^i(pl_0) = pl_i$ and so on. The areas enclosed by the segments of the stable and the unstable manifolds connecting two successive PIPs are called *lobes*. Denote the lobes enclosed by segments of the stable and the unstable manifolds connecting $pl_i, ql_i$ by $Dl_i$, and the ones that between $ql_i, pl_{i+1}$, by $El_i$, when again, $Dl_i = F^i(Dl_0), El_i = F^i(El_0)$ (the equivalent notation is used for the right side), see figure 6.

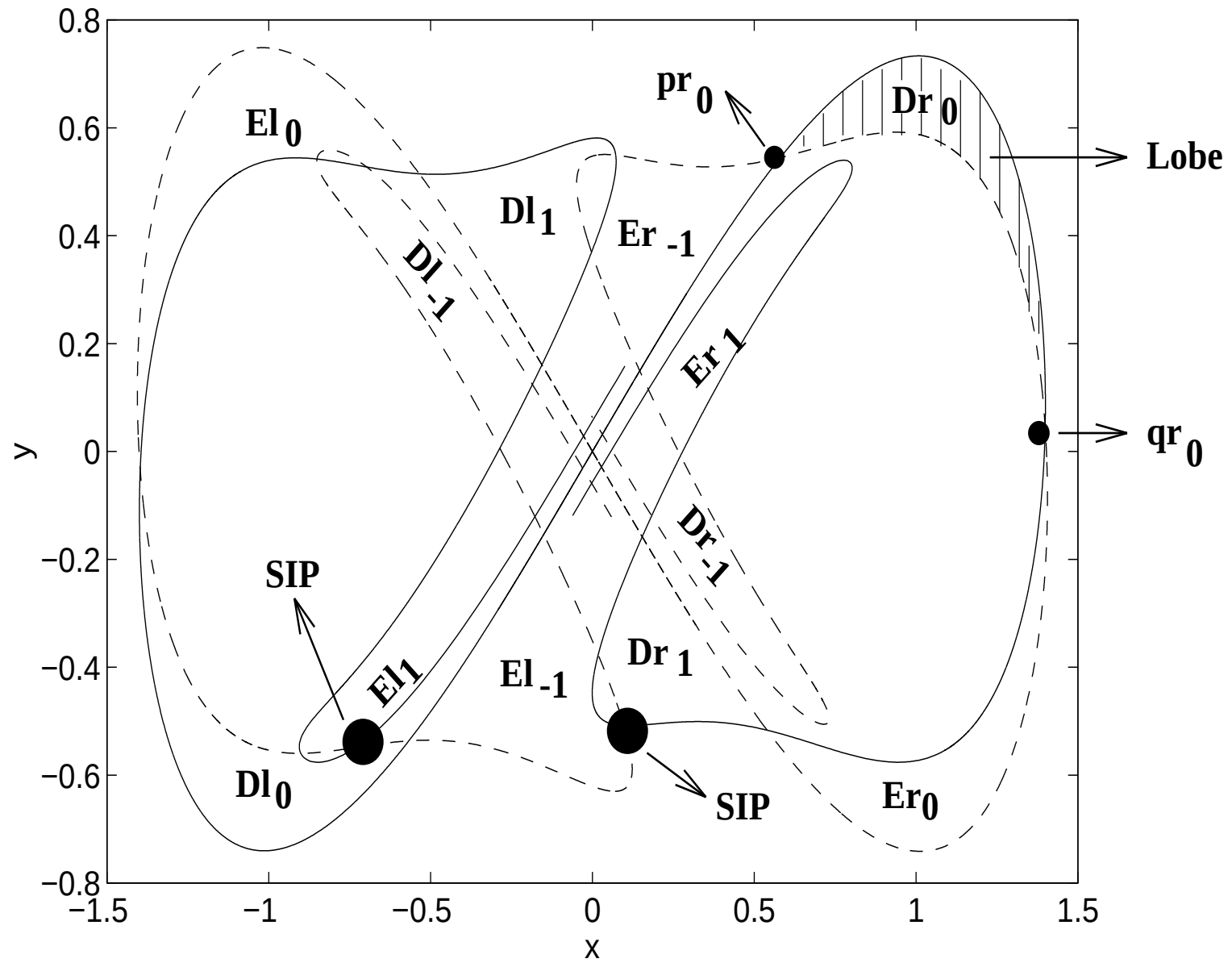

- -: stable manifold, —: unstable manifold

Figure 6: Illustration of Lobes, PIPs and SIPs.

If $El_j \cap Dl_0 \neq \emptyset$ (respectively $Er_j \cap Dr_0 \neq \emptyset$) for some non negative integer $j$, or if $Dl_{k+1} \cap Er_0 \neq \emptyset$ (respectively $Dr_{k+1} \cap El_0 \neq \emptyset$) for some non negative integer $k$, then there exist Secondary Intersection Points (SIPs) in these intersections. The integers $j, k$ are the *transition numbers* of the corresponding SIPs [8]. The minimal transition numbers (the minimal integers $j, k$) for which this happens on the left side (respectively right side) of the hyperbolic fixed point are called the *the structural indices* $\ell_{ll}, \ell_{lr}$ (respectively $\ell_{rr}, \ell_{rl}$) of the homoclinic tangle [8]. Namely, these structural indices are exactly the *transition numbers* of the secondary homoclinic points which belong to the first intersection of the corresponding lobes. For example, in figure 6 the structural indices are: $\ell_{ll} = 1$, $\ell_{lr} = 1$, $\ell_{rl} = 1$ and $\ell_{rr} = 2$. Each such structural index imposes minimal complexity for the structure of the homoclinic tangle. For example, the length growth rate of line segments along the unstable manifold (the topological entropy) increases as the indices decrease. Thus, dividing the parameter space to regions according to the values of these indices corresponds, approximately, to a plot of "level sets" of the topological entropy. The dividing curves correspond to secondary homoclinic bifurcations. Large SA which are not associated with primary homoclinic tangencies seems to appear only in regions in which at least one structural index is less or equal to 1 (see also [21]).

The perturbational method for calculating the secondary homoclinic bifurcation curves is described below. For simplicity it is presented specifically for the AFDO model. More generally, it may be applied to nearly Hamiltonian dissipative systems, which satisfy some generic assumptions (see [20]).

Consider the Secondary Melnikov Function - SMF (see [31, 32, 33] and [20]):

$$h_2^{cd}(t_0, \varepsilon) = M_c(t_0) + M_d(t_{1cd}(t_0, \varepsilon)); \quad c, d \in \{l, r\}, \tag{17}$$

where $M_c(t)$ is the Melnikov function, and $t_{1cd}$ is defined by:

$$t_{1cd}(t_0, \varepsilon) = \begin{cases} t_0 + P(\varepsilon M_c(t_0)), & M_c(t_0) < 0 \\ t_0 + \frac{1}{2} P(\varepsilon M_c(t_0)), & M_c(t_0) > 0 \end{cases} ; c, d \in \{l, r\}, \tag{18}$$

$P(H)$ is the period of the unperturbed periodic orbit with energy $H$, and $H = 0$ on the separatrix. For sufficiently small $\varepsilon$, simple zeros (respectively degenerate zeros) of (17) imply transverse secondary homoclinic intersections (respectively tangencies) with a *transition number*:

$$j_{cd}(t_0, \varepsilon) = \left[ \frac{t_{1cd}(t_0, \varepsilon)}{T} \right] - s(t_0); \quad 0 \leq t_0 < T, \tag{19}$$

where $[x]$ is the integer part of $x$, $T = \frac{2\pi}{\omega}$ is the period of the perturbation and $s(t_0)$ is either 0 or 1, depending on the interval to which $t_0$ belongs (equation (24) below). The structural index $\ell_{cd}$ $(c, d \in \{l, r\})$ is defined to be the minimal transition number $j_{cd}(t_0, \varepsilon)$. For sufficiently small $\varepsilon$ this analytical definition of the structural index meets the geometrical definition described above [33]. It follows that typically [5] a change in the structural index may be found at a bifurcation point for equation (17). Indeed, under some generic conditions on $h_2^{cd}(t_0, \varepsilon; \mu, \delta)$, the structural indices satisfy [33]:

$$\ell_{cd} = j_{cd}(t_{0cd}, \varepsilon_{cd}); \quad c, d \in \{l, r\}, \tag{20}$$

where $(t_{0cd}, \varepsilon_{cd})$ ($\varepsilon_{cd}$ small) are the solutions to the standard equations for a bifurcation point:

$$h_2^{cd}(t_0, \varepsilon) \;=\; 0; \quad 0 \le t_0 < T \tag{21}$$

$$\frac{\partial h_2^{cd}(t_0, \varepsilon)}{\partial t_0} \;=\; 0, \tag{22}$$

defined in the appropriate time-interval for $t_1$:

$$[\ell_{cd} + s(t_{0cd})]T \;\le\; t_{1cd}(t_{0cd}, \varepsilon_{cd}) < [\ell_{cd} + s(t_{0cd}) + 1]T; \tag{23}$$

$$s(t_{0cd}) \;=\; \begin{cases} 0; & t_{0cd} \in \left[0, \frac{T}{2}\right) \\ 1; & t_{0cd} \in \left[\frac{T}{2}, T\right) \end{cases}, \tag{24}$$

where $t_{1cd}(t_0, \varepsilon)$ is defined by equation (18).

Typically, for $\varepsilon_{cd}$ sufficiently small, one finds a sequence of two bifurcation values, $\varepsilon_{cd}^1 < \varepsilon_{cd}^2$. The corresponding solutions $(t_{0cd}^i, \varepsilon_{cd}^i)$; $i = 1, 2$ of (21) and (22) divide the parameter space to three regions: below the hyper-surface $\varepsilon = \varepsilon_{cd}^1(\omega, \gamma, \beta, \delta, \ell_{cd})$ there are no SIPs, between the hyper-surfaces $\varepsilon = \varepsilon_{cd}^1(\omega, \gamma, \beta, \delta, \ell_{cd})$ and $\varepsilon = \varepsilon_{cd}^2(\omega, \gamma, \beta, \delta, \ell_{cd})$ two SIPs occur (see for example the intersections denoted by an arrow in figure 7), and above the hyper-surface $\varepsilon = \varepsilon_{cd}^2(\omega, \gamma, \beta, \delta, \ell_{cd})$ two additional SIPs occur (see for example the intersections of the lobes above the origin in figure 7).

Moreover, equations (17) - (22) may be brought to a simple form, as shown in appendix A. Using the asymptotic expansion for the period function $P(H)$ (equation (6)), these equations may be solved analytically if both $\beta = 0$ and

[5] Another source for changes in the structural indices are points of discontinuity of $t_1(t_0)$, see [20].

$\delta = 0$, to find approximations to the secondary homoclinic bifurcation points $(t_0^i, \varepsilon_{cd}^i(t_0^i))$, $c, d \in \{l, r\}$, $i \in \{1, 2\}$. The solutions for $\beta = 0$, $\delta = 0$ may be used to solve these equations for $\beta \neq 0$ or/and $\delta \neq 0$ small by the use of asymptotic expansions in powers of $\beta$ and $\delta$ for $i = 1$, and in powers of $\beta^{\frac{1}{3}}$ and $\delta^{\frac{1}{3}}$ for $i = 2$. To get more accurate results, using $P(H)$ instead of its approximate form, we use a Newton method, combined with a linear prolongation scheme (see appendix A and [20] for details). Plotting the bifurcation values, $\varepsilon_{cd}^i(\omega, \gamma, \beta, \delta, \ell_{cd})$; $c, d \in \{l, r\}$, $i \in \{1, 2\}$, obtained from the solutions to these equations, for $\gamma, \beta, \delta$ fixed and varying $\omega$, gives the secondary homoclinic bifurcation curves in parameter space $(\omega, \varepsilon)$, labeled by the structural indices, $\ell_{cd} = 0, 1, 2, ..., m < \infty$, as in figure 7a.

A simple lower bound to the homoclinic bifurcation curves (compare with (30)) is given by:

$$\bar{\varepsilon}_{cd}(\omega, \gamma, \beta, \delta, \ell_{cd}) = \begin{cases} \dfrac{P^{-1}(\frac{2\pi}{\omega}(\ell_{cd}+1))}{\max_{t_0}|M_c(t_0)|}, & c = d \\ \\ \dfrac{P^{-1}(\frac{4\pi}{\omega}(\ell_{cd}+1))}{\max_{t_0}|M_c(t_0)|}, & c \neq d \end{cases} \tag{25}$$

Therefore, using (6,7), we find

$$\varepsilon_{cd}^{1,2}(\omega, \gamma, \beta, \delta, \ell_{cd}) \geq \frac{16 e^{-\frac{2\pi}{\omega}(\ell_{cd}+1)}}{\max\{|\gamma F_c - \frac{4}{3}\delta|, |\gamma F_c + \frac{4}{3}\delta|\}}. \tag{26}$$

These give a simple lower bounds on the secondary homoclinic bifurcation curves[6]. Moreover, geometrically these lower bounds correspond to the values of $\varepsilon$ for which the lobes may get involved in a $1 : (\ell_{cd} + 1)$ resonance (see below).

### 2.3.1 Comparison between numerical and analytical results

The analytical method described above for finding the secondary homoclinic bifurcations is of a perturbational nature. Thus, as proved in [33], in the

[6]Note that since the approximation to leading order in $H$ for the period function, $P(H)$, is used here to calculate $\bar{\varepsilon}_{cd}$, we get that $\bar{\varepsilon}_{cd} \equiv \bar{\varepsilon}_{cc}$ (see the appendix for more details). Hence, the curve $\bar{\varepsilon}_{cc}(\ell_{cc} = n)$ serves as a lower bound to all the eight secondary homoclinic bifurcation curves, related to the structural index $\ell_{cc} = n$ , $\varepsilon_{cd}^{1,2}(\ell_{cc} = n)$, with $c, d \in \{l, r\}$ and $n \geq 0$.

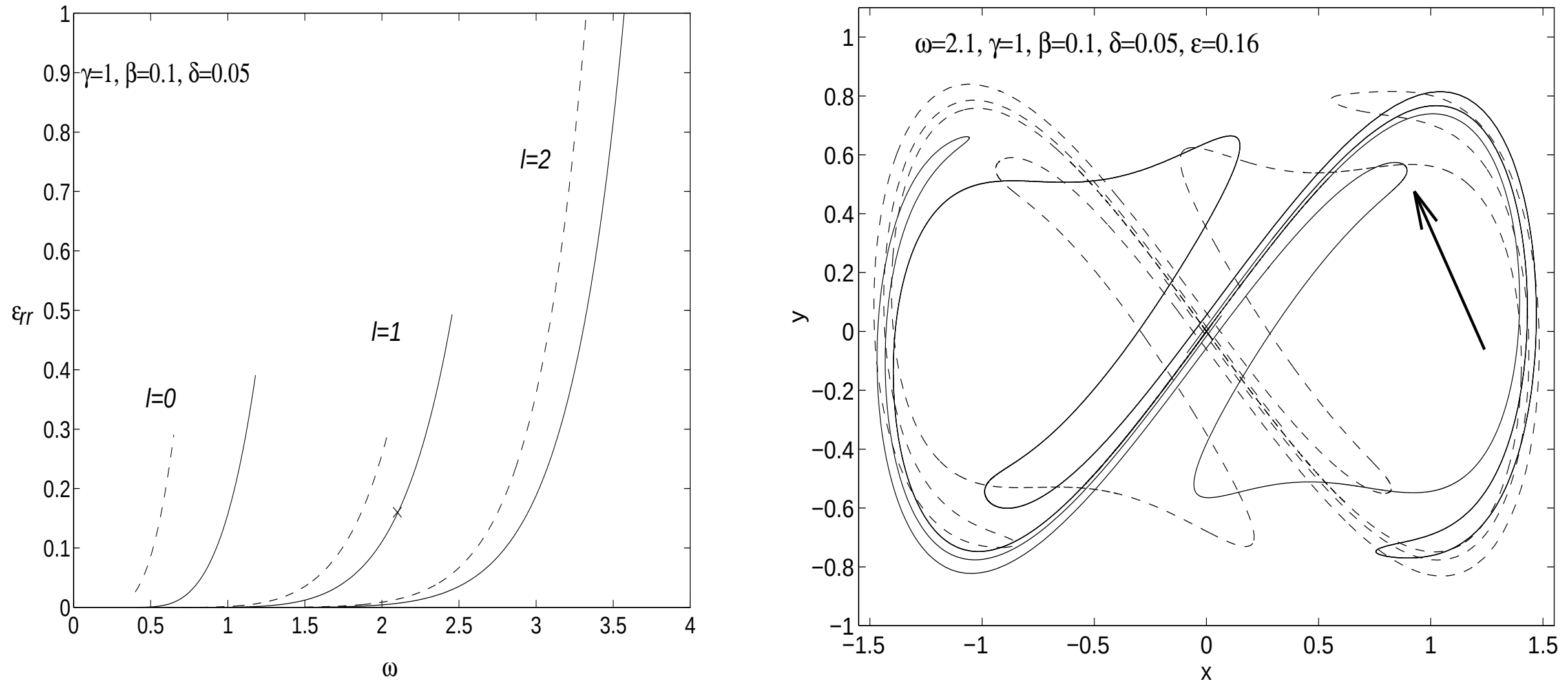

**a** **b**

—: $\varepsilon^1_{rr}$, - -: $\varepsilon^2_{rr}$, - -: stable manifold, —: unstable manifold.

Figure 7: Secondary homoclinic bifurcation curves and the homoclinic tangle.

limit of small $\varepsilon$ values it is guaranteed to supply a good approximation to the actual bifurcation value. Here we examine how good of an approximation the analytical formulae supply for finite $\varepsilon$ values. Indeed, excellent agreement is achieved between the analytical predictions for the occurrence of SIPs and the numerical calculations of the stable and unstable manifolds for $\varepsilon$ as large as 0.3 and $\ell_{cd} \geq 1$; $c, d \in \{l, r\}$, see for example figure 7. In this figure the $'\times'$ at the $(\omega, \varepsilon)$ parameter space indicates the parameter values for which the manifolds, presented on the right figure, are calculated. On the right figure, the corresponding near-tangency of the manifolds is indicated by an arrow. In fact, the larger the $\ell$'s, the larger the $\varepsilon$ values for which the zeroth-order approximation is found to be adequate. For example, for $\ell = 2$ we find excellent agreement up to $\varepsilon \approx 1$. This is not surprising since large $\ell$'s (and finite $\varepsilon$) correspond to large $\omega$'s for which the Melnikov function coefficient becomes exponentially small, thus the effective perturbation is small.

For $\ell_{cc} = 0$ the agreement between the numerical and the analytical results is not as favorable (notice that this is a finite $\varepsilon$ effect: letting $\varepsilon \to 0$, with all other parameters held fixed, necessarily implies that $\ell \to \infty$); This is due to the passage of the manifolds through a $1:1$ resonance relation between the periodic orbits inside the homoclinic loop and the forcing pe-riod $2\pi/\omega$; Namely, the energy level to which the manifolds are pushed by the Melnikov function is near the energy level for which a $1:1$ resonance occurs[7] Now, the construction of the SMF uses the Whisker map (see [6, 8, 10, 30, 31, 32, 33, 40]) in which the motion of the interior orbits is approximated by unperturbed periodic motion. This approximation fails near a $1:1$ resonance. Hence, the analytical approximation for the $\ell_{cc} = 0$ bifurcation curve are inaccurate even for small values of $\varepsilon$. Notice that in this limit $\omega$ is varied with $\varepsilon$, hence this observation is not contradictory to the SMF theorems which hold in the limit $\varepsilon \to 0$ with all other parameters held fixed [33]. Indeed, to avoid passage of the manifolds through a $1:1$ resonance, $\varepsilon$ should satisfy the condition: $\varepsilon_{cc} < \frac{P^{-1}(\frac{2\pi}{\omega})}{M_c(t_{0cc})}$. This condition holds for $\ell_{cc} \geq 1$. For the outer indices $\ell_{cd}, c \neq d$ the problem of $1:1$ resonance was not encountered.

# 3 Strange attractors.

In this section numerical evidence for the existence of SA, and observations regarding their location in parameter space and their structural properties in phase space are presented.

## 3.1 Numerical scheme for detecting SA.

Simple numerical experiments showing Poincaré maps of the AFDO (using DSTOOLS [7]) suggested that SA appear in the area of the parameter space related to the structural indices $\ell_{cd} = 0$; $c, d \in \{l, r\}$. To investigate this subject more thoroughly the Lyapunov exponents of orbits of (3) were computed (see [38]). Viewing (3) as an autonomous system, each orbit has three Lyapunov exponents, one zero, one negative and the third may be either positive or negative. A positive third Lyapunov exponent indicates the existence of a strange attractor (see [9]), while a negative third Lyapunov exponent indicates that the orbit is attracted to a periodic sink.

An efficient stopping criteria for the Lyapunov exponents calculation is developed, using the distinction between SA or sinks with long transients and simple sinks. First, to remove transient behavior, $N_{in}$ iterations of the Poincaré map[8] are calculated. Then, every $N_{it}$ iterations ($N_{it} = 100$) a line is fitted to the logarithm of the modulus of the last $N_{it}$ Lyapunov exponent values. The program stops if one of the following events occurs:

1. The largest Lyapunov exponent is negative. Then there exist a periodic sink. The exact value of the negative Lyapunov exponent is not sought.

2. The largest Lyapunov exponent is positive and the slope of the fitted line is nearly zero (up to an error of $1e - 6$).

3. The total number of Poincaré map iterations exceeds 10000. In this case no decision is made regarding the existence or non-existence of an attractor. In practice the stopping criteria 1 and 2 occur before 10000 iterates are computed.

[7]Indeed, the $1:m$ resonance relation for the periodic orbits of equation (3) is given by: $P(H) = \frac{2\pi m}{\omega}$ . Since by definition of $t_{1cc}$ (see equation (18)), $P(\varepsilon_{cc} M_c(t_{0cc})) = t_{1cc} - t_{0cc}$, and by condition (23), $t_{1cc} - t_{0cc} \in [\ell_{cc}\frac{2\pi}{\omega}, (\ell_{cc}+1)\frac{2\pi}{\omega})$, the manifolds "pass" through the $1:1$ resonance zone for $\ell_{cc} = 0$.

In case 2, when a positive Lyapunov exponent is detected, all the Lyapunov exponents are calculated. Hence, the Lyapunov dimension [24] of the strange attractor may be calculated. Note that the Lyapunov dimension, $D_L$ , is an upper limit for the Capacity (or Box-Counting) dimension, $D_0$ [24].

In the numerical experiments, $\gamma, \beta, \delta$ are fixed and $\varepsilon$ and $\omega$ are varied along and near the secondary homoclinic bifurcation curves $\varepsilon^i_{cd}(\omega, \gamma, \beta, \delta, \ell_{cd})$; $\quad i = 1, 2; \quad c, d \in \{l, r\}$, of section 2.3.

## 3.2 Windows of SA.

For various parameter values, numerical evidence suggests the existence of SA in "windows" in the parameter space. These windows are aligned near the secondary homoclinic bifurcation curves which are related to the structural indices $\ell_{ll}, \ell_{lr} = 0$, see figure[9] 8, and the magnification of the windows in figure 9. While theoretically such regions should appear near all tangent bifurcation with arbitrary $\ell$, we did not detect in our numerical search any SA near the bifurcation curves with $\ell > 0$. This suggests that the size of the parameter regions for which SA appear decreases dramatically with $\ell$.

[8] $N_{in} = 200$ was found sufficient.

[9] Notice that in figure 8a $\beta = 0$ hence $\varepsilon^i_{ll} \equiv \varepsilon^i_{rr}$ and $\varepsilon^i_{lr} \equiv \varepsilon^i_{rl}$, whereas in figure 8b $\varepsilon^i_{ll} \neq \varepsilon^i_{rr}$ and $\varepsilon^i_{lr} \neq \varepsilon^i_{rl}$, but, for clarity, $\varepsilon^i_{rr}$, $\varepsilon^i_{rl}$ are not plotted. In figure 8c $\varepsilon^i_{rr}$, $\varepsilon^i_{rl}$ are not defined for the specified $\omega$ values (since the Melnikov function $M_r(t_0)$ has no zeros).

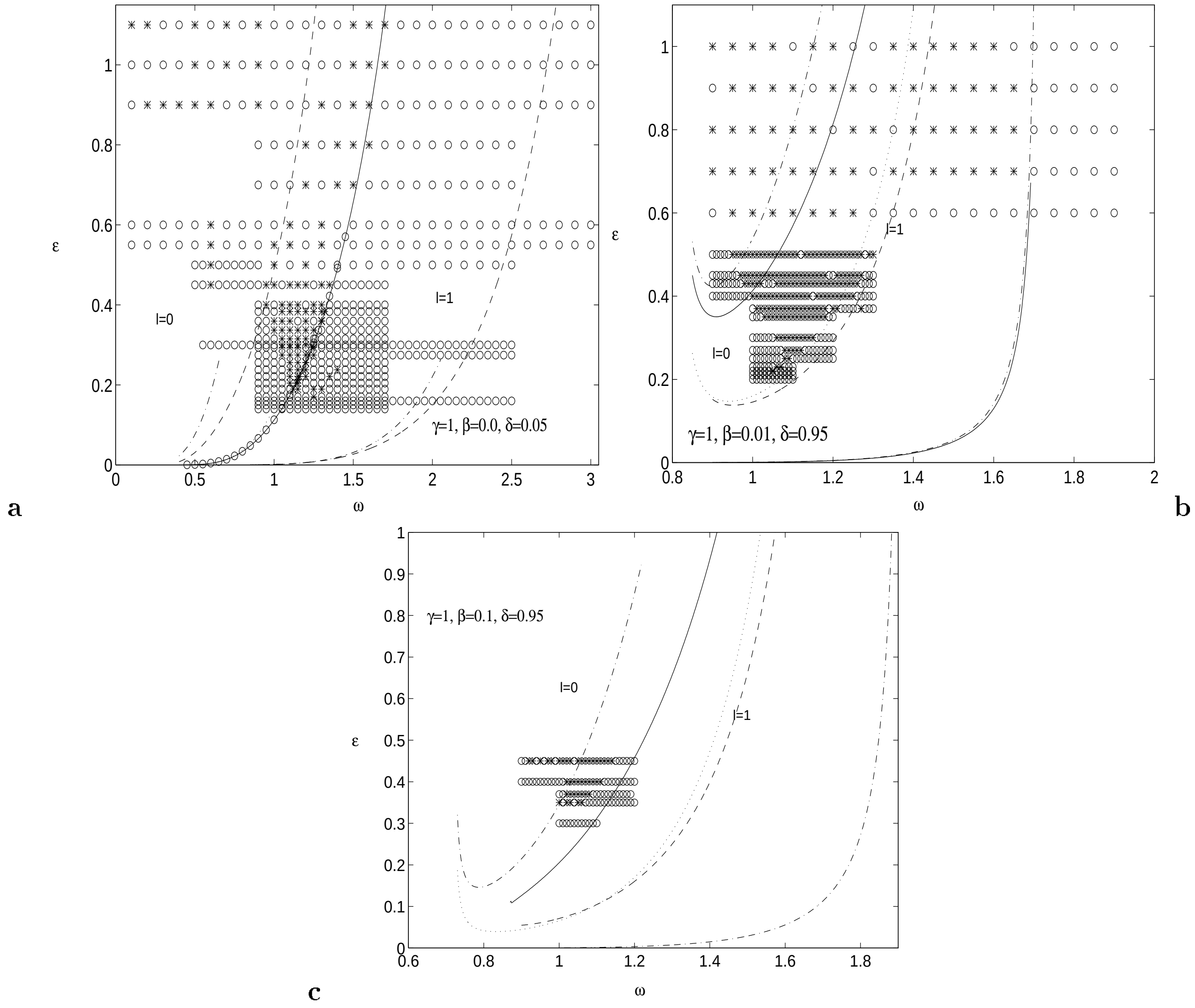

Figure 8: Secondary homoclinic bifurcation curves and SA for the AFDO. $'*'$ indicates a strange attractor (positive Lyapunov exponent), 'o' indicates a periodic orbit (negative Lyapunov exponent), —: $\varepsilon^1_{lr}$, - -: $\varepsilon^2_{lr}$, $\cdots$: $\varepsilon^1_{ll}$, $-\cdot-$: $\varepsilon^2_{ll}$.
a) $\beta = 0$, $\delta = 0.05$.
b) $\beta = 0.01$, $\delta = 0.95$.
c) $\beta = 0.1$, $\delta = 0.95$.

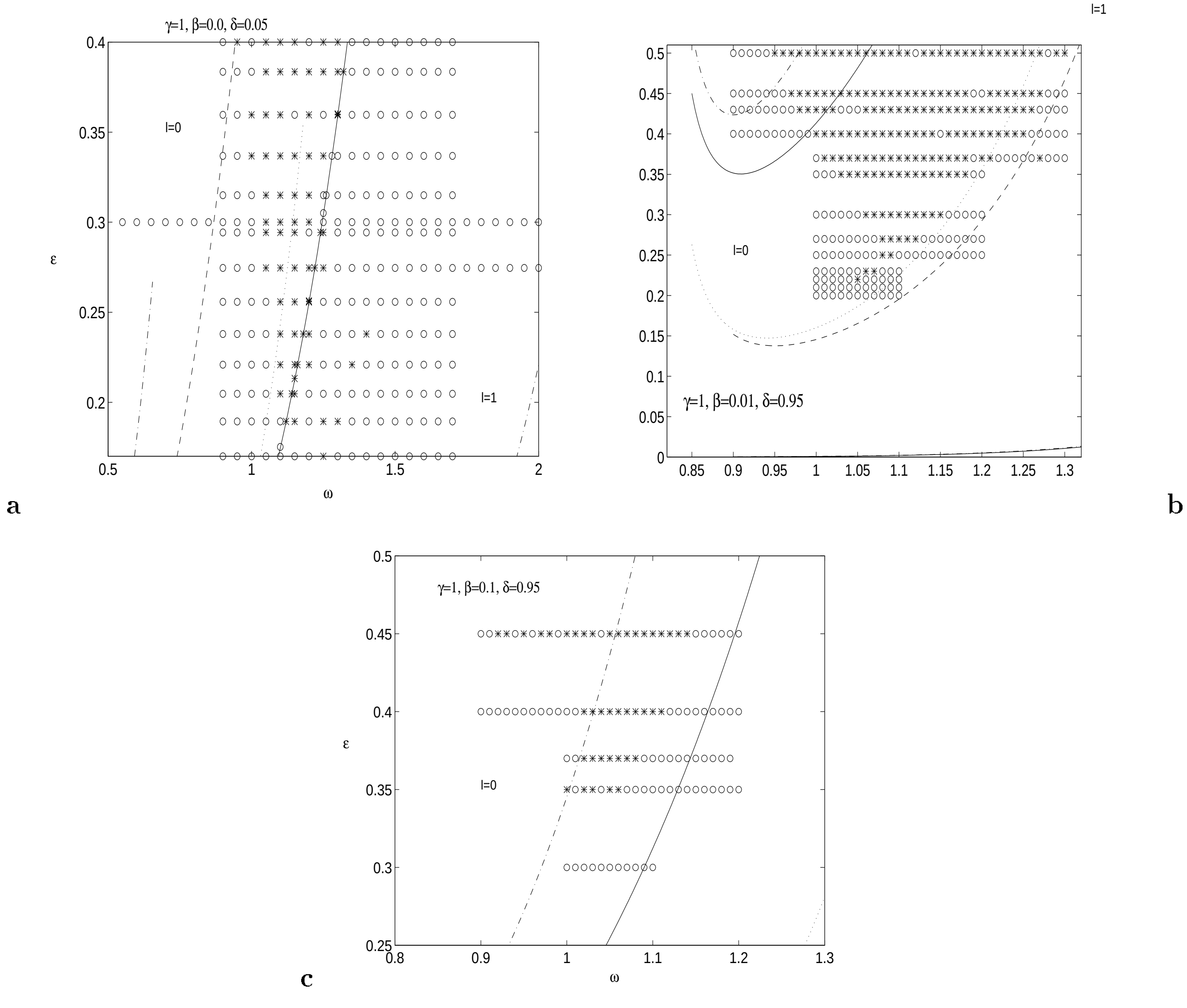

Figure 9: Secondary homoclinic bifurcation curves and SA for the AFDO. Magnification of certain regions of figures a-c.

In fact, our perturbational methods for detecting the homoclinic tangencies associated with the structural indices $\ell_{ll} = 0$ are inaccurate, see section 2.3. Nonetheless, the predicted analytical bifurcation curves for $\ell_{cd} = 0; \quad c, d \in \{l, r\}$ still lie in the area of parameter space near which the actual bifurcation curves exist. Moreover, observe that the SA appear only in the region of the parameter space $(\omega, \varepsilon)$, which is above the second secondary bifurcation curve $\varepsilon^2_{ll}(\omega, \gamma, \beta, \delta, \ell_{ll} = 1)$ (see section 2.3). For $\varepsilon < 1$, our predictions for this curve are accurate, hence, this curve may be considered as a lower bound to the region in parameter space in which SA appear, see figures 8.

Another feature of the SA windows is that they all seem to appear above a threshold value $\tilde{\varepsilon}(\omega, \gamma, \beta, \delta) \geq 0.2$. Namely, they do not seem to extend to the small $\varepsilon$ values to which some of the $\ell = 0$ bifurcation curves extend. Numerical calculations of the stable and the unstable manifolds of the origin for the minimal $\varepsilon$ values for which strange attractors are found suggest that this curve is a specific homoclinic bifurcation curve; above this curve the lobe $El_1$ intersects the lobe $Dl_0$ at five or more (six, seven or eight) homoclinic intersection points. See for example figure 10. Notice that the homoclinic bifurcation curve $\varepsilon^2_{ll}(\omega, \gamma, \beta, \delta, \ell_{ll} = 1)$ in the $(\omega, \varepsilon)$ parameter space gives a lower bound to this curve.

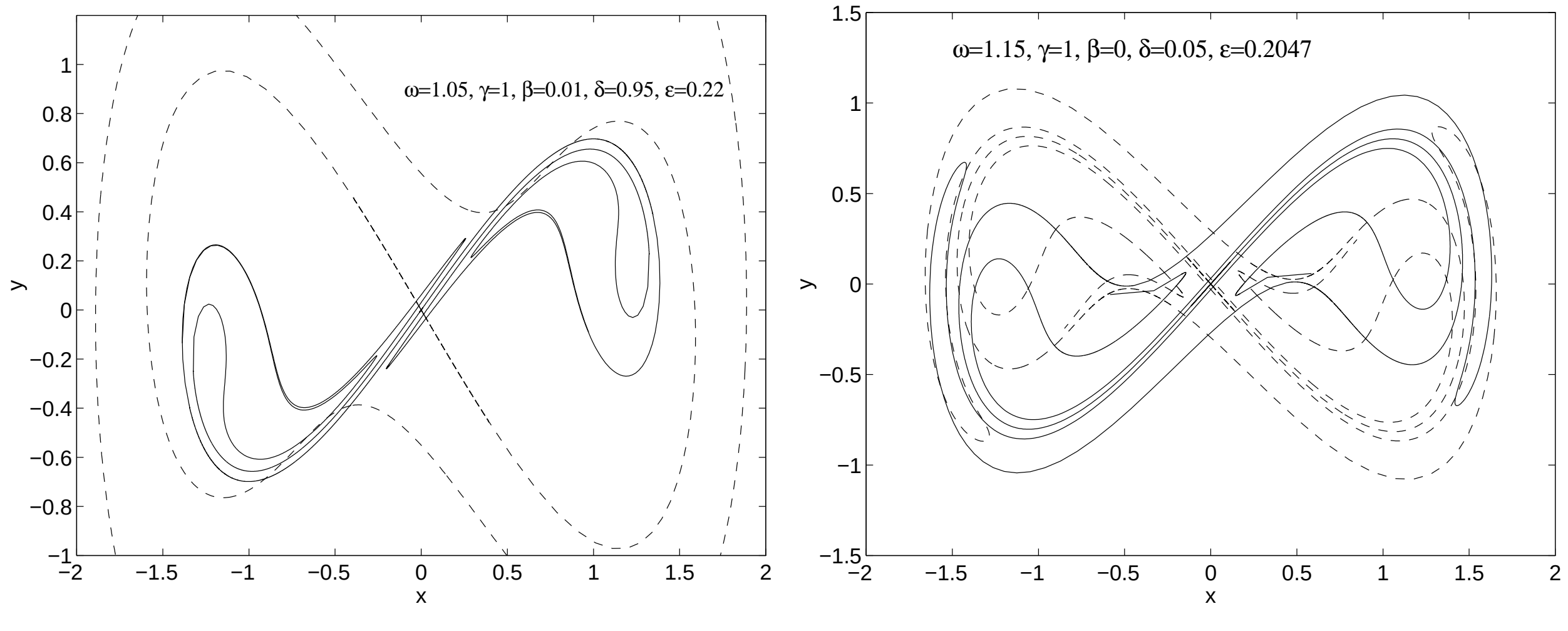

Figure 10: The stable and the unstable manifolds corresponding to SA.

The structure of the SA that are obtained vary with the parameters. The main forms of the attractors which were found are described next.

## 3.3 The phase space structure of the SA.

The observed SA have the following two distinct properties:

1. The attractors may be one sided (i.e. the attractor is contained in the right/left half phase space plane) *or* two sided (with one or two components).
2. The attractors may have strong dissipative features *or* may have nearly conservative features.

The first property depends on the location of the parameter values with respect to the division to regions I, II and III. One sided SA may appear in region II near secondary homoclinic bifurcations or near the boarders between the regions, namely near primary homoclinic tangencies, see figure 11. There, the transition between two-sided SA, denoted by 'T', and one sided SA, denoted by 'O', is shown[10]. In fact, near the boarder between regions II and III three different SA may appear: one sided SA on the right half plane, *two* one-sided SA or a two-sided SA. Near the boarder between regions I and II left sided SA, coexisting with a sink on the right half plane, were observed.

The second property seems to depend mainly on the ratio $\delta/\gamma$ and is roughly independent of the other parameters (in regions where SA exist); This is somewhat surprising since the area contraction per Poincaré map is given by $\exp(-\varepsilon\delta 2\pi/\omega)$ - thus strong dissipation may be achieved for fixed $\delta/\gamma$ by increase of $\varepsilon$, without essential changes in the structure of the SA. For $\frac{\delta}{\gamma} \ll 1$, the two sided strange attractors seems to have nearly conservative features of a chaotic region, see figure 12. These features persist in the window shown in figure 9a, and even when $\varepsilon = 1$, though the attractor is more structured, it has "fat" regions in which no filamentation is observed. The positive Lyapunov exponent corresponding to this figure is $\log_2(\alpha_1) \approx 0.1987$, and the corresponding Lyapunov dimension is $D_L \approx 1.9036$. The SA indicated in figure 8a are of such structure. Most of these attractors have a positive Lyapunov exponent of about 0.2. The maximal observed deviation from this value is 0.02. Nearly conservative one-sided SA were not observed.

[10]This transition is not continuous in $\frac{\delta}{\gamma}$; Between the value of $\frac{\delta}{\gamma}$ for which a two sided SA appears, and the value of $\frac{\delta}{\gamma}$ for which a one sided SA appears, there may be some values of $\frac{\delta}{\gamma}$ for which no SA appear.

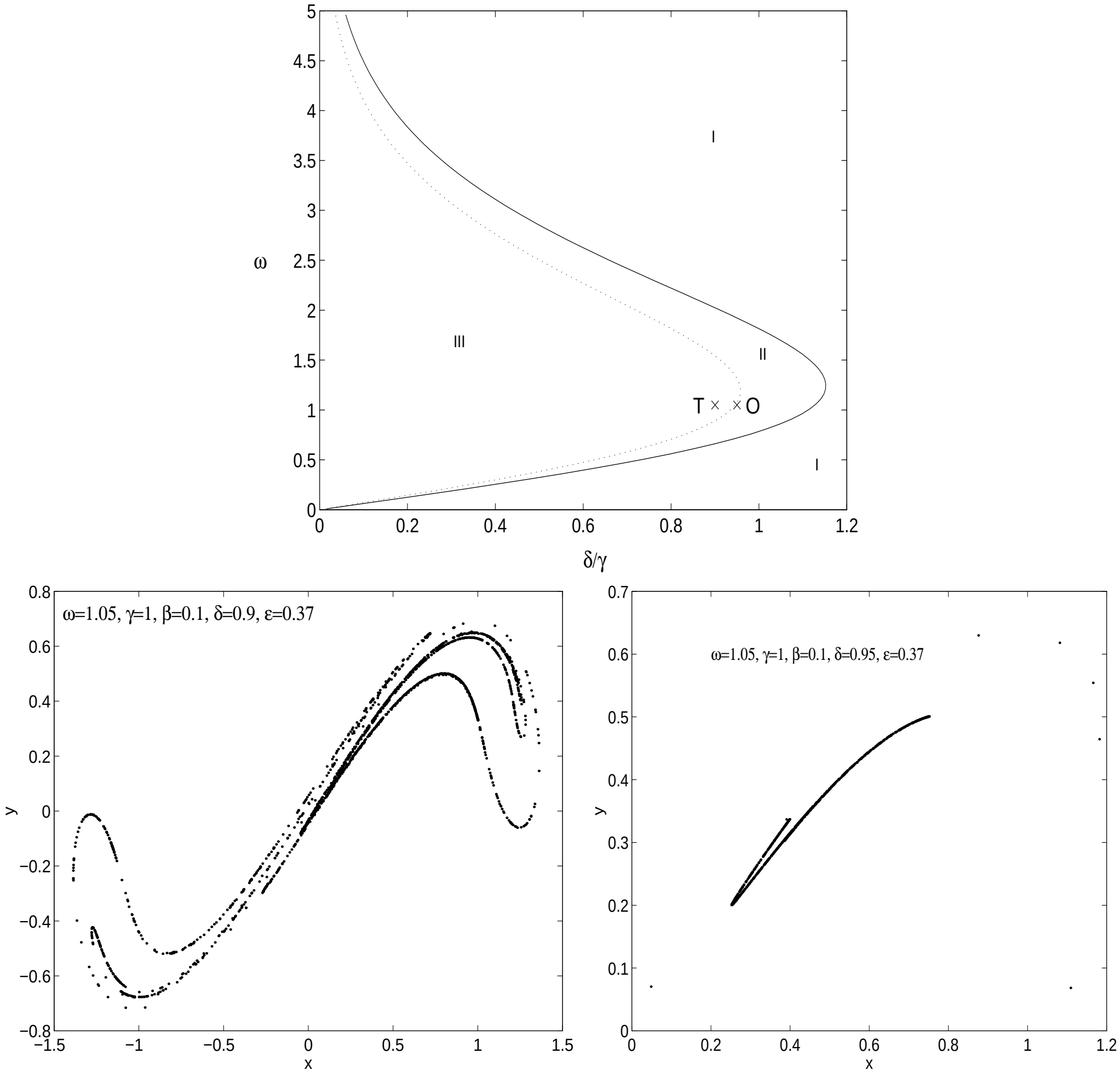

T O

Figure 11: Transition from two sided to one sided strange attractor.

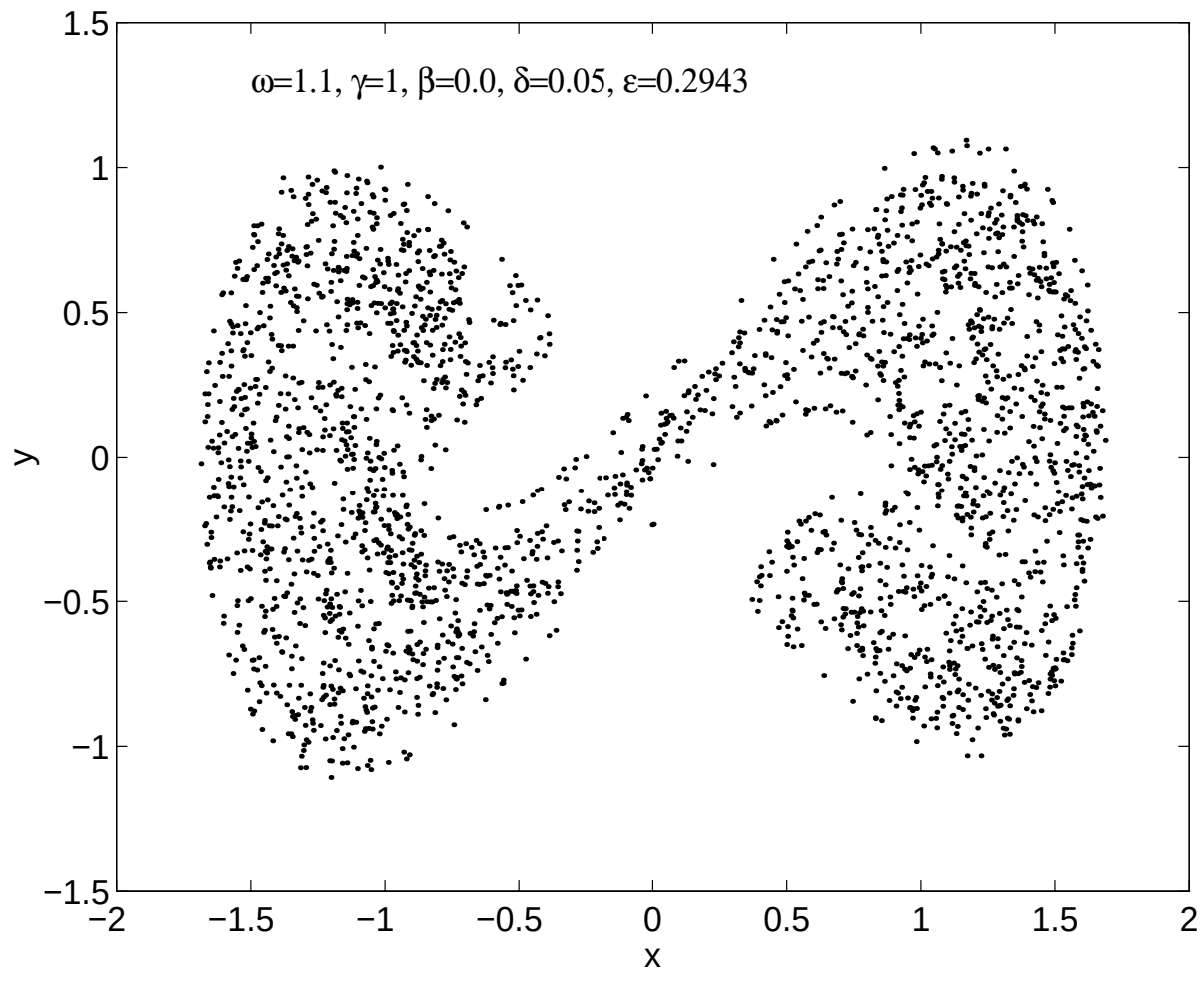

Figure 12: Nearly-conservative strange attractor.

When $\frac{\delta}{\gamma}$ is not small ($\delta/\gamma = 0.2$ is already in this region), the structure of the SA is associated with the folding of the unstable manifold as in strongly dissipative systems, see figure 13. The corresponding positive Lyapunov exponent is $\log_2(\alpha_1) \approx 0.1741$, and the Lyapunov dimension is $D_L \approx 1.2975$. Comparing with the above results for the nearly-conservative attractors, we observe that the Lyapunov exponent is less sensitive to the attractor's structure than the Lyapunov dimension. The SA indicated in figure 9b are of strongly dissipative nature; some are two-sided with a positive Lyapunov exponent very close to 0.17, with maximal deviation of 0.03, and some are one-sided. In figure 14 such one-sided SA is shown, its positive Lyapunov exponent is $\log_2(\alpha_1) \approx 0.0862$, and its Lyapunov dimension is $D_L \approx 1.1453$. The SA presented in figure 8c are all one-sided strongly dissipative (Hénon-like) SA. r The values of the positive Lyapunov exponent are $0.08 \pm 0.04$, about half of the Lyapunov exponents of the two-sided SA.

The Lyapunov exponent and dimension of the attractors seems to be quite robust. The dependence of the Lyapunov exponents on the values of $\varepsilon$ along a secondary homoclinic bifurcation curve related to a structural index $\ell_{cd} = 0$ is shown in figure 15. An example for the dependence of the Lyapunov dimensions on the values of $\varepsilon$ along such bifurcation curve is shown in figure 16. The plunges in the figure correspond to parameter values for which no SA exist.

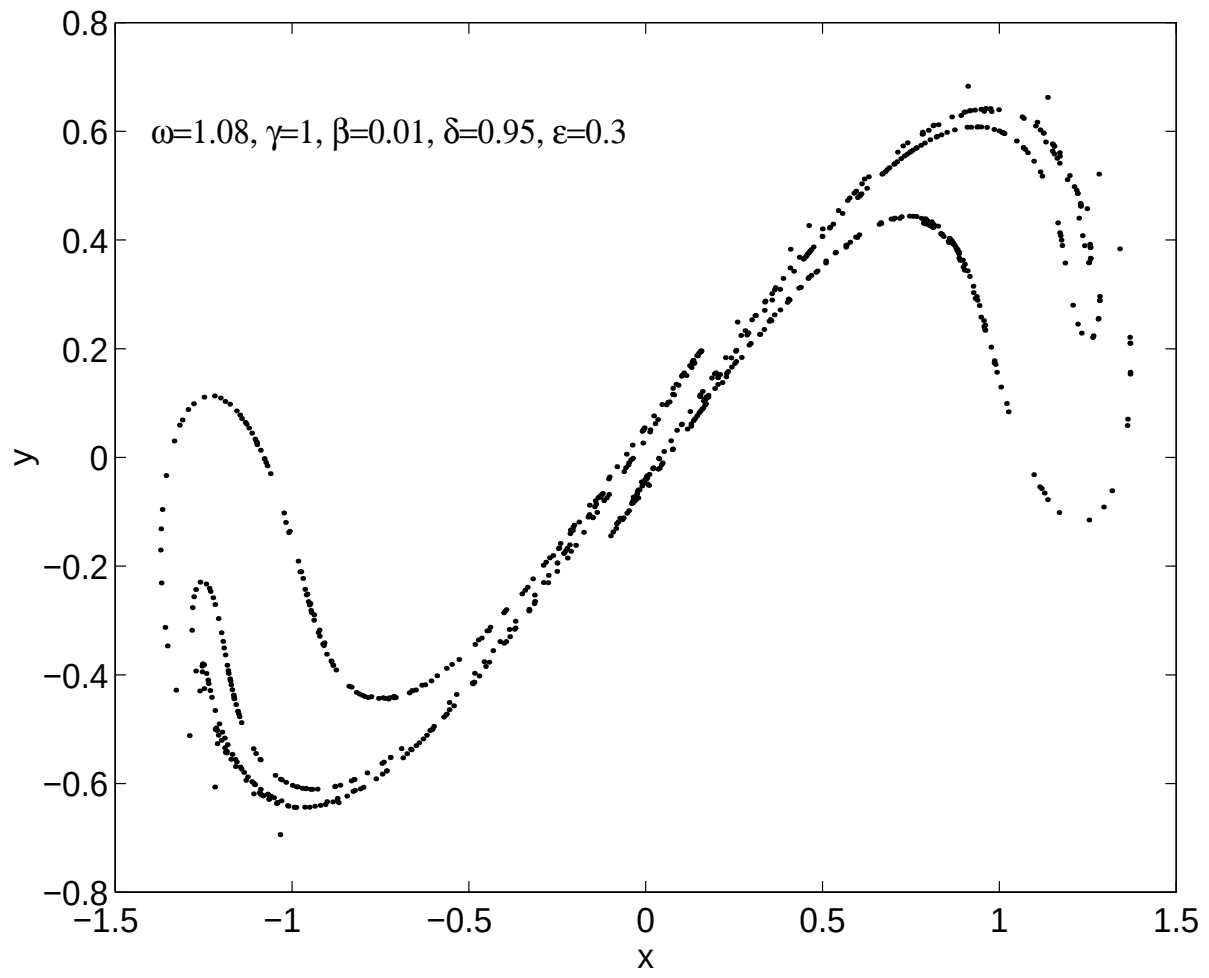

Figure 13: Two sided strange attractor of the AFDO, with structure referring to large values of the relation $\frac{\delta}{\gamma}$.

In [20] symbolic dynamics of segments of the unstable manifold is constructed for the AFDO, and for general dissipative systems which unfold homoclinic tangencies. From this symbolic dynamics a transfer matrix may be constructed for each set of the structural indices $\ell_{cd}$; $c, d \in \{l, r\}$. It follows from [28] and [32], that $\log(\lambda)$, where $\lambda$ is the modulus of the largest eigenvalue of the transfer matrix, gives a lower bound on the topological entropy of the Poincaré map. The lower bound on the topological entropy for the AFDO with $\ell_{cd} = 0$; $c, d \in \{l, r\}$, corresponding to the region in parameter space where two sided SA appear, is $\log_2(3.9231) = 1.9720$. The lower bound on the topological entropy for the AFDO, corresponding to the existence of one sided SA is $\log_2(3.6709) = 1.8761$.

The above results are consistent with the inequalities describing the relations between topological entropy, entropy (Kolmogorov-Sinai invariant) and positive Lyapunov exponents (see [9] and [39]):

$$h(\rho) \leq h_{top}$$
$$h(\rho) \leq \sum_{\lambda_i > 0} \lambda_i$$

Where, $\rho$ is an ergodic measure with compact support, with respect to a diffeomorphic map $F$, and $\lambda_i$ are the positive Lyapunov exponents corre-

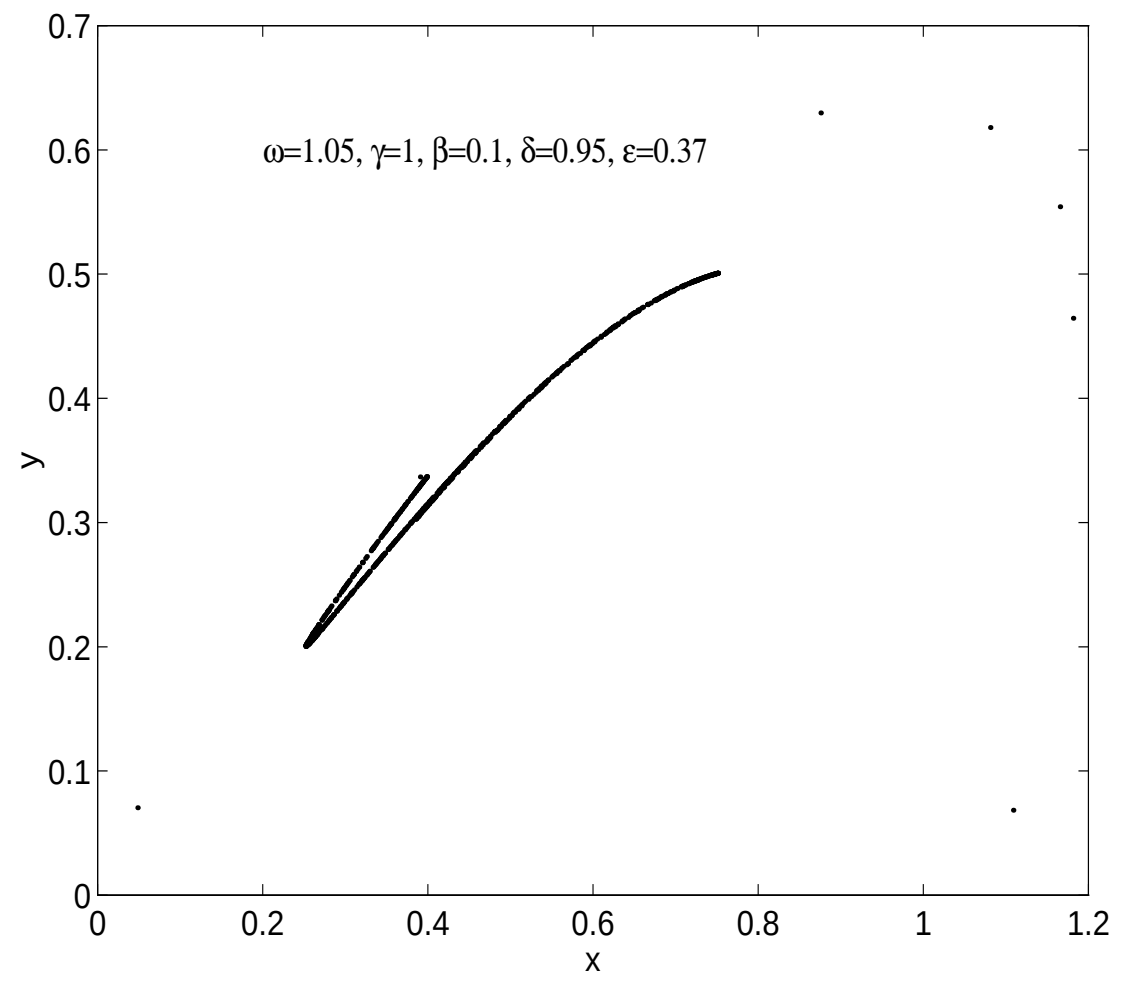

Figure 14: One sided strange attractor.

sponding to a dense orbit of $F$. i.e. for the structural indices $\ell_{cd} = 0$ we get:

$$h(\rho) \leq \lambda_1 \approx \log_2(\alpha_1) < \log_2(\lambda) \leq h_{top}.$$

Notice that the lower bounds obtained for the topological entropy are by order of magnitude larger then the corresponding Lyapunov exponents and that there is a very little difference between the one-sided and the two-sided cases. See [31] and [20] for the construction of the symbolic dynamics of the lobes, and details on how the transfer matrices and the lower bounds on the topological entropy may be calculated.

## 3.4 Basins of attraction.

Initial conditions may be attracted to the various attractors which exist in the phase space. Here, we distinguish between three types of attractors: those located entirely on the left (respectively right) half plane and those which are located on both sides of the $y$ axis. We do not distinguish here between the basins of attraction of different sinks or SA, see [26, 36, 35, 15, 21] for detailed study of these issues.

The flux of phase space area into the left/right sides may be calculated

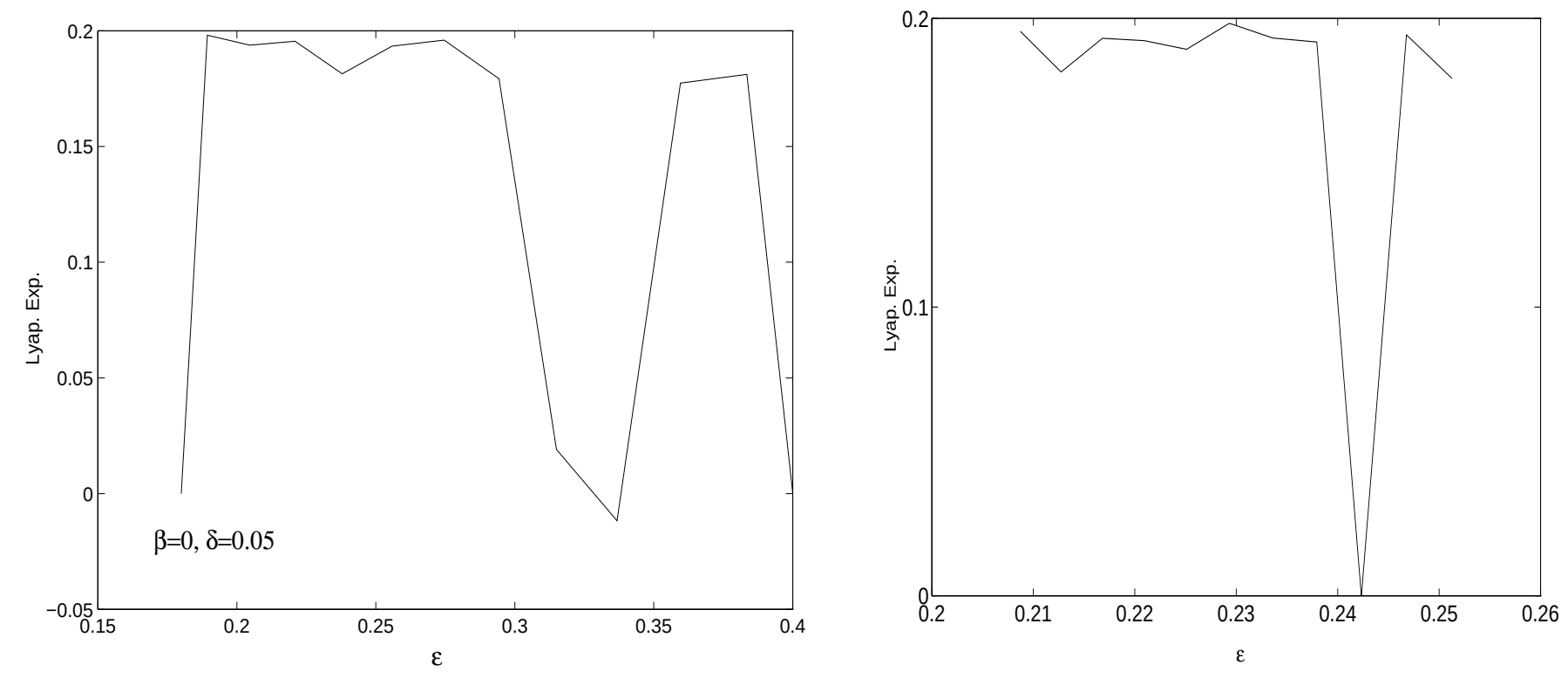

**a** **b**

Figure 15: The positive Lyapunov exponent variation.
a) Samples of $\varepsilon$ values along the bifurcation curve $\varepsilon = \varepsilon^1_{lr}(\omega, \gamma, \beta, \delta; \ell_{lr})$ with $\beta = 0, \quad \frac{\delta}{\gamma} = 0.05$ and $\ell_{lr} = 0$.
b) Zoom in on $\varepsilon$ values in the interval $[0.2, 0.25]$.

to first order in $\varepsilon$ by integrating the Melnikov function:

$$
\begin{aligned}
\Delta^{in}_{l,r} &= \int_0^{s_0} M_{l,r}(s)ds + \int_{s_1}^{T} M_{l,r}(s)ds = \qquad (27)\\
&= \begin{cases} -\frac{8\pi\delta}{3\omega}; \quad s_0 = s_1, \;\; M_{l,r}(s) < 0 \quad \forall s \in [0, T) \\ \\ -\frac{2\gamma}{\omega} F_{l,r}\sqrt{1 - \left(\frac{4\delta}{3\gamma F_{l,r}}\right)^2} - \frac{4\delta}{3\omega}(2\arcsin(\frac{4\delta}{3\gamma F_{l,r}}) + \pi); \quad s_0 \neq s_1, \end{cases}\\
s_0 &= \frac{1}{\omega}\arcsin(\frac{4\delta}{3\gamma F_{l,r}}), s_1 = \frac{\pi}{\omega} - s_0,
\end{aligned}
$$

where $s_0, s_1$ are determined by $M_{l,r}(s) \leq 0$ for $s_0 \leq s \leq s_1 \in [0, T]$. Thus, on the $n^{th}$ iterate the initial phase space area, $\left|\Delta^{in}_{l,r}\right| \exp(2\pi\varepsilon\delta n/\omega)$ , is swept into the left/right side. In region I, $M_{l,r}(s) < 0 \quad \forall s \in [0, T)$, and $|\Delta^{in}_r| = |\Delta^{in}_l|$ to order $\varepsilon$. Since $F_r < F_l$ it can be easily shown that $|\Delta^{in}_r| > |\Delta^{in}_l|$ in regions II and III. Thus, for parameter values corresponding to these regions, the influx to the right side is always larger than the influx to the left side (recall that $\beta > 0$).

If $M_{l,r}(s)$ has simple zeroes (so $s_0 \neq s_1$) then, similarly, the flux out of the left/right sides is given by:

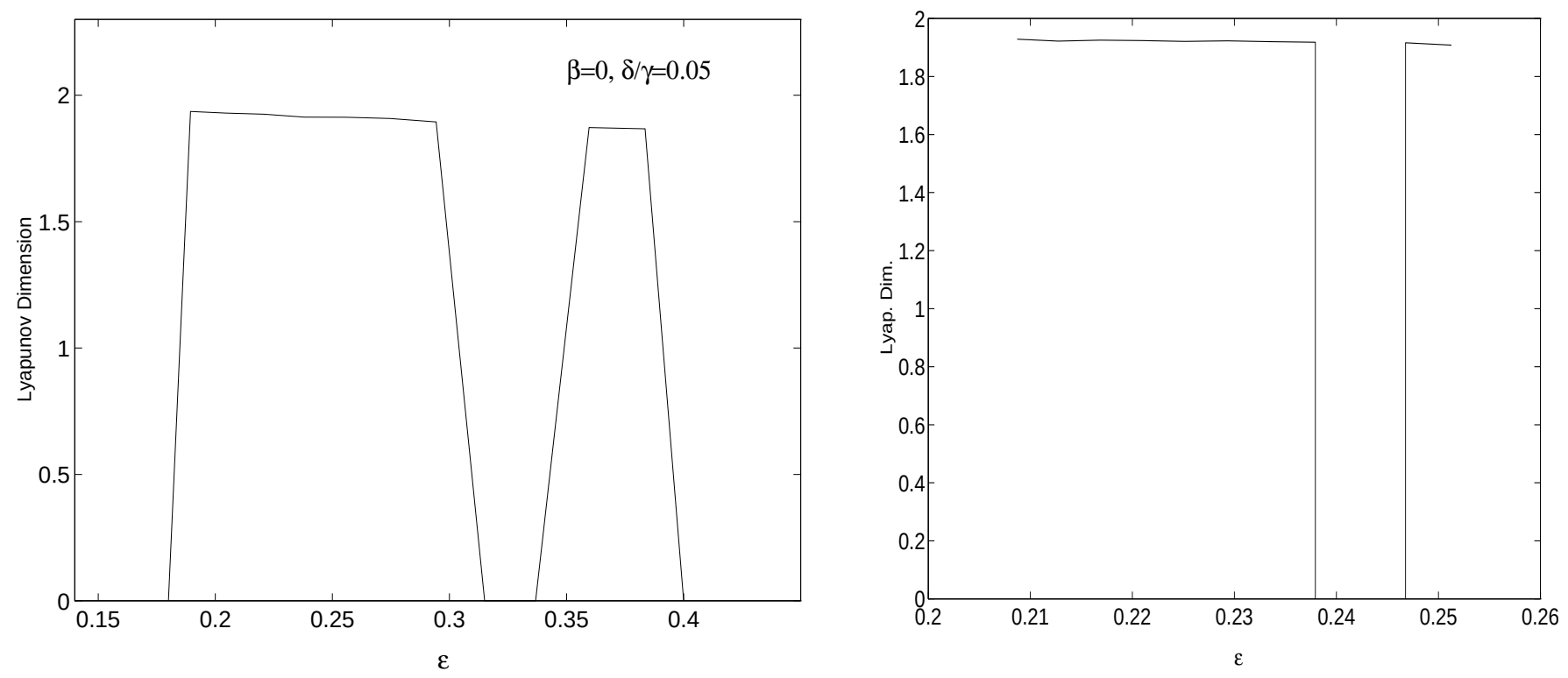

**a** **b**

Figure 16: The Lyapunov dimension variation.
a) Samples of $\varepsilon$ values along the bifurcation curve $\varepsilon = \varepsilon^1_{lr}(\omega, \gamma, \beta, \delta; \ell_{lr})$ with $\beta = 0$, $\frac{\delta}{\gamma} = 0.05$ and $\ell_{lr} = 0$.
b) Zoom in on $\varepsilon$ values in the interval $[0.2, 0.25]$.

$$\begin{aligned} \Delta^{out}_{l,r} &= \int_{s_0}^{s_1} M_{l,r}(s)ds = \\ &= \frac{2}{\omega}(\gamma F_{l,r}\sqrt{1-\left(\frac{4\delta}{3\gamma F_{l,r}}\right)^2} + \frac{4\delta}{3}\arcsin(\frac{4\delta}{3\gamma F_{l,r}}) - \frac{2\pi\delta}{3}). \end{aligned} \tag{28}$$

In region II, $\Delta^{out}_r = 0$, hence obviously $|\Delta^{out}_l| \geq |\Delta^{out}_r|$ in this region, and it can be shown that $|\Delta^{out}_l| \geq |\Delta^{out}_r|$ in region III as well (i.e. for $\delta \leq \frac{3\gamma F_r}{4}$).

One might expect that the ratio between the fluxes to the right and left regions determines the ratio between the sizes of the basins of attraction. However, in regions II and III near the boarder line between the regions, this picture may change dramatically; There are cases for which all the initial conditions which are numerically integrated are attracted to the right side.

Numerical calculations of the basins of attraction of the left/right attractors suggest a more detailed description:

1. No intersections on both sides - region I:

   To order $\varepsilon$, by equation (27) the left and right influxes areas are equal, but high order terms alter these results. Indeed, numerically it is found that the right basin is larger than the left one. Moreover, it's area seems to grow monotonically as $\delta$ decreases to it's threshold value, $3\gamma F_l/4$.

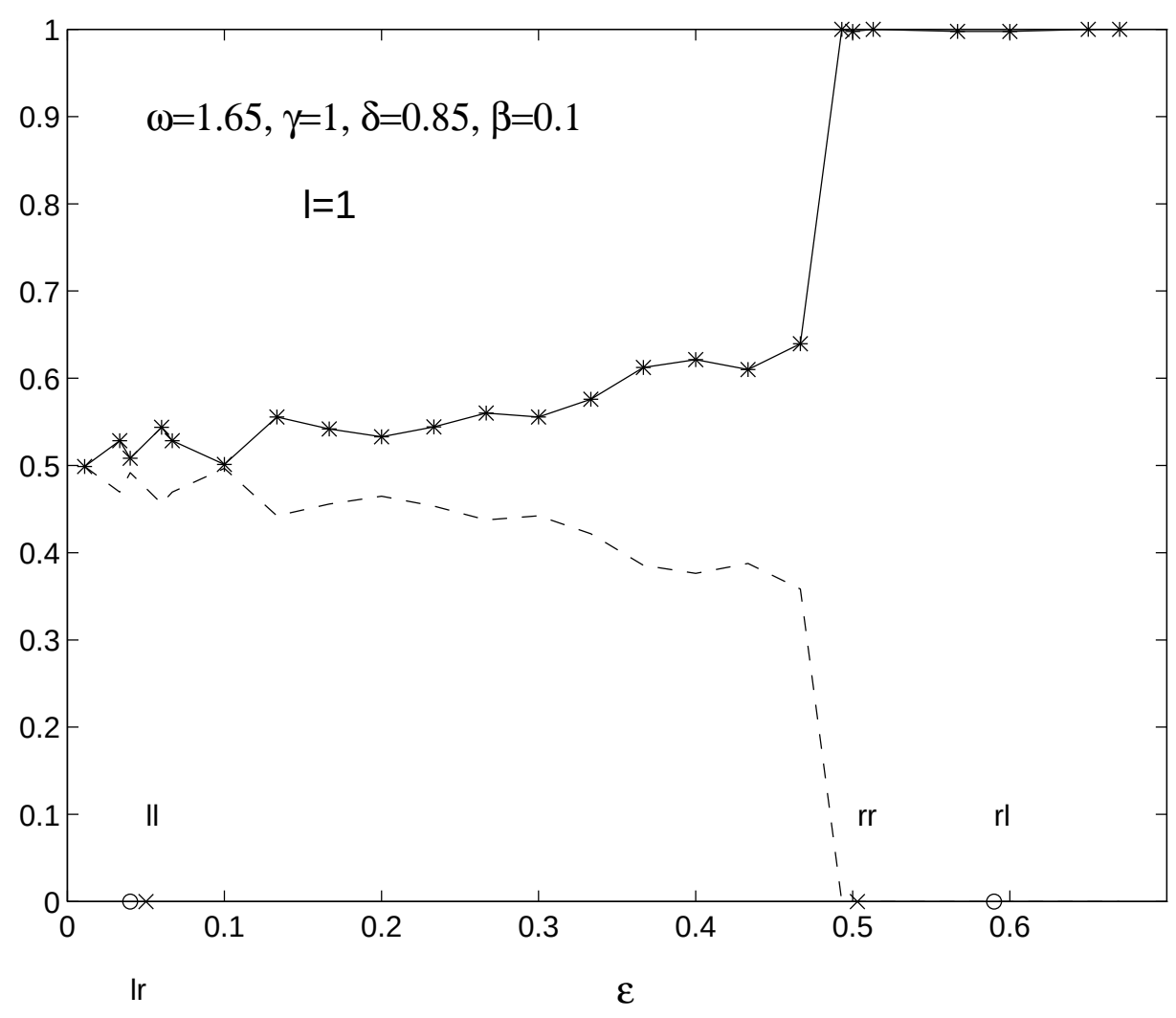

Figure 17: The relative size of the right and left basins of attractions.
—: the right basin size, - - -: the left basin size, ∗: sampled values of $\varepsilon$. The secondary homoclinic bifurcation values $\varepsilon(\ell_{cc}=1)$ are denoted by 'x', and the secondary homoclinic bifurcation values $\varepsilon(\ell_{cd}=1)$ are denoted by 'o'.

2. Intersections on the left side or on both sides - regions II and III:

   As $\delta$ is further decreased, the area of the right basin continues to grow. However, the growth in the basin's area seems to be discontinuous. This phenomenon is associated with the "boundary metamorphosis" [15] of the subharmonics. For some parameter values in region III a *two sided* periodic orbit or SA may exist, hence in these cases, some (or most) of the initial conditions are not attracted neither to a left nor to a right attractor.

There seems to be a correlation between the discontinuities and the occurrence of homoclinic tangencies. See for example figure 17. In this figure the percentage of the sampled phase space area which is swept to the left/right sides, up to an approximated error of $\pm 0.03$, is presented as a function of $\varepsilon$ (all the other parameters are fixed), were the values of $\varepsilon$ which correspond to secondary homoclinic bifurcation curves are specified. The results of this figure are somewhat puzzling; One would expect that left/right attractors are composed, for $\ell \geq 0$, from the attracting resonances. Hence, when the unstable manifold intersects through the resonance region, it depletes the resonance, thus decreasing the relative area of the basin of attraction of the corresponding side. However, in figure 17, along with results that confirm this scenario (see for example the jump near $\varepsilon = 0.06$), we observe quite the contrary results (see for example the huge jump near $\varepsilon = 0.5$). Possibly other, undetected bifurcation is responsible to these results.

In [21], numerical results regarding the relation between homoclinic and other bifurcation curves and the basins of attraction of systems with a cubic like potential well are presented. There, it has been suggested that the bifurcation curve, corresponding to what we call here a secondary homoclinic bifurcation curve with a structural index $\ell = 0$, is of great significance, since closely beneath it they numerically observed a chaotic escape (i.e. a destruction of the basin of attraction of the SA). Their chaotic escape corresponds, in a case of a *closed*[11] system, to a decrease in the basin of attraction of one side and an increase in the basin of attraction of the other side. Thus, the current results are in agreement with the results obtained in [21] for *open* systems. Possibly, the critical curves (near which SA appear or loose stability) that they have observed numerically, correspond to the curves discussed above: the homoclinic bifurcation curve above which the lobe $El_1$ intersects the lobe $Dl_0$ at at-least five points, and the curve at which the unstable manifold intersects the resonance region.

# 4 Discussion and conclusions

The qualitative differences between the flows under symmetric and asymmetric forcing loom when primary homoclinic intersections/tangencies occur only on one side of the saddle fixed point (region II + its neighborhood). This region may be of significant size even for very small asymmetry values ( $\beta \ll 1$) if the forcing frequency is appropriately set. It is of negligible size in the adiabatic limit, hence, to the best of our knowledge, was not observed in previous works which have considered asymmetric potentials with adiabatic forcing. In this fat region II, the system may posses one one-sided-SA (strange attractors), two one-sided-SA, or one two-sided SA. In the former case the attractor may be situated on either side of the fixed point; however,

[11] A system is called *closed* if some forward iteration of the Poincaré map of a segment of the unstable manifold which has left the left/right region returns to it. See [20] for more precise definition.

the nature of the basin of attraction of the left and right SA seems to be different.

We find that the relative size of the basin of attraction to the left/right attractors is usually not sensitive to its strangeness (i.e. the size of the basin does not change significantly when a SA is destroyed/created). In general, the basin of attraction of the right attractor is always a bit larger than that of the left attractor, where in most cases the difference between the fluxes to the right and left regions determines the ratio between the sizes of the basins of attraction. This occurs in a continuous and natural way in region I, However, near the border between regions II and III, where one sided SA may appear (on either side of it), this picture may change dramatically; There are cases for which all the initial conditions which are numerically integrated are attracted to the right side.

The robust, observable (hence physically significant) SA appear near primary homoclinic tangencies and near secondary homoclinic tangencies with small structural index. A key perturbational tool for finding the latter is the SMF (secondary Melnikov function [33]). Generally, it is found that the SMF supplies excellent analytical prediction to the occurrence of secondary homoclinic tangencies even for relatively large values of $\varepsilon$. However, it fails near the $\ell = 0$ homoclinic bifurcations, exactly in the region were robust SA exist. Thus, only lower bounds and approximate curves for the regions were SA are observed are found. We suspect that both phenomena (the failure of the SMF and the appearance of SA) are associated with the involvement of the stable and unstable manifolds in a $1:1$ resonance. Thus we derived a simple lower bound for $\varepsilon$ above which the manifolds enter the $1:1$ resonance was derived (eq. (25)). The study of the relation between the resonance, the manifolds and the strange attractor, and the construction of a more accurate approximation to the homoclinic bifurcation curves near a $1:1$ resonance are left for future work.

The structure of the SA vary with the parameters; as the ratio $\frac{\delta}{\gamma}$ increases the values of the positive Lyapunov exponent slightly decreases and the Lyapunov (fractal) dimension decreases significantly. Surprisingly, we find that fixing this ratio and varying the other parameters in one of the "windows" for which SA exist, the structure, Lyapunov exponent and Lyapunov dimension of the SA hardly change. Such a variation does change, in particular, the dissipation (area contracting) rate per Poincaré map.

## Acknowledgments

We thank A. A. Vasiliev and Dmitry Turaev for discussion and comments. This research is supported by MINERVA Foundation, Munich/Germany.

# A Finding secondary homoclinic bifurcation curves.

In this appendix some technical aspects regarding the method of solution of ( 17) - (21) are described.

From (17) - (21) we construct the equations:

$$t^i_{1cd}(t_0) = M_d^{-1,i}(-M_c(t_0)) + j_{cd}T, \tag{29}$$

$$\varepsilon^i_{cd}(t_0) = \begin{cases} \dfrac{P^{-1}(t^i_{1cd}(t_0) - t_0)}{M_c(t_0)}, & c = d \\ \dfrac{P^{-1}(2(t^i_{1cd}(t_0) - t_0))}{M_c(t_0)}, & c \neq d \end{cases} \tag{30}$$

where $c, d \in \{l, r\}$; $i = 1, 2$; $j_{cd} \in \mathbf{N}$, and from equation (7) for the Melnikov function of the AFDO, one gets:

$$M_d^{-1,1}(x) = \frac{1}{\omega}\arcsin\left(\frac{x}{\gamma F_d(\omega,\beta)} + \frac{4\delta}{3\gamma F_d(\omega,\beta)}\right),$$

$$M_d^{-1,2}(x) = \frac{\pi}{\omega} - \frac{1}{\omega}\arcsin\left(\frac{x}{\gamma F_d(\omega,\beta)} + \frac{4\delta}{3\gamma F_d(\omega,\beta)}\right).$$

Notice that $M_d^{-1,i}(x)$ are undefined for $\left|\frac{x}{\gamma F_d(\omega,\beta)} + \frac{4\delta}{3\gamma F_d(\omega,\beta)}\right| > 1$ , hence there are some values of $t_0$ in $[0, T)$ for which $\varepsilon^i_{cd}(t_0)$ from (30) are undefined. From (20) and (23) we get that for $t_0 \in [0, T)$: $j_{cd} = \ell_{cd} + s(t_0)$, where $s(t_0) = 0$ for $t_0 \in [0, \frac{T}{2})$, and $s(t_0) = 1$, for $t_0 \in [\frac{T}{2}, T)$. The above equations ((29) and (30)) to be dependent on the perturbation parameters: $\beta$ (asymmetry), $\delta$ (dissipation) $\gamma$ (the amplitude of the forcing) and $\omega$ (the frequency of the for forcing). In addition, for $|H| = |\varepsilon M_c(t_0)| \ll 1$ we get (see equations (6), (21), (22) and (29)):

$$0 = \frac{\partial h_2^{cd}(t_0, \varepsilon)}{\partial t_0} \approx M_c'(t_0) + M_d'(M_d^{-1,i}(-M_c(t_0)))\left(1 - \frac{M_c'(t_0)}{M_c(t_0)}\right); \tag{31}$$

$c, d \quad \in \quad \{l, r\}, \quad i = 1, 2.$

Hence, for the AFDO, equation (29) becomes:

$$t_{1cd}^1(t_0) = \frac{1}{\omega}[2\pi + \arcsin(-\frac{F_c(\omega,\beta)}{F_d(\omega,\beta)}\sin(\omega t_0) + \frac{8\delta}{3\gamma F_d(\omega,\beta)})] \tag{32}$$

$$+(\ell_{cd} + s(t_0))\frac{2\pi}{\omega},$$

$$t_{1cd}^2(t_0) = \frac{1}{\omega}[\pi - \arcsin(-\frac{F_c(\omega,\beta)}{F_d(\omega,\beta)}\sin(\omega t_0) + \frac{8\delta}{3\gamma F_d(\omega,\beta)})] \tag{33}$$

$$+(\ell_{cd} + s(t_0))\frac{2\pi}{\omega},$$

$$s(t_0) = \begin{cases} 0, & t_0 \in [0, \frac{\pi}{\omega}) \\ 1, & t_0 \in [\frac{\pi}{\omega}, \frac{2\pi}{\omega}) \end{cases} ; \quad c, d \in \{l, r\},$$

and equation (30) becomes:

$$\varepsilon_{cd}^i(t_0) \approx \begin{cases} \dfrac{-16\exp(t_0 - t_{1cd}^i(t_0))}{\gamma F_c(\omega)\sin(\omega t_0) - \frac{4\delta}{3}}; & c = d, \\ \\ \dfrac{16\exp(t_0 - t_{1cd}^i(t_0))}{\gamma F_c(\omega)\sin(\omega t_0) - \frac{4\delta}{3}}; & c \neq d, \end{cases}$$

$$i = 1, 2. \tag{34}$$

These approximations are valid for sufficiently small $H's$ of the period function (see equation (6)), namely for:

$$\varepsilon_{cc}^i(t_0)[\gamma F_c(\omega)\sin(\omega t_0) - \frac{4\delta}{3}] \rightarrow 0- \tag{35}$$

$$\tag{36}$$

$$\varepsilon_{cd}^i(t_0)[\gamma F_c(\omega)\sin(\omega t_0) - \frac{4\delta}{3}] \rightarrow 0+$$

$$\tag{37}$$

And, for such sufficiently small values of $H$, equation (31) is:

$$\frac{F_c(\omega,\beta)}{F_d(\omega,\beta)}\sin(\omega t_0)+\cos(\omega t^i_{1cd}(t_0))[\tan(\omega t_0)$$
$$-\frac{\omega\sin(\omega t_0)}{\sin(\omega t_0)-\dfrac{4\delta}{3\gamma F_c(\omega,\beta)}}]\approx 0, \tag{38}$$

where $F_c, F_d$ are as in (9), (10) of section 2.1, for $c,d\in\{l,r\}$, and $\varepsilon^i_{cd}(t_0)$ is calculated here with the use of the approximated value of the period function of AFDO, $P(H)$, from (6).

**Remark 1** *Since the exact inverse function of the period function of the AFDO, $P^{-1}(x)$, cannot be found analytically, $(t^i_0, \varepsilon^i_{cd}(t^i_0))$; $i=1,2$ are found by solving the equations:*

$$P(\varepsilon M_c(t_0)) = \tau^i_{cd}(t_0) \tag{39}$$

$$\frac{\partial P(\varepsilon M_c(t_0))}{\partial t_0} = \frac{\partial \tau^i_{cd}(t_0)}{\partial t_0}, \tag{40}$$

*where:*

$$\tau^i_{cd}(t_0) = \begin{cases} t^i_{1cc}(t_0)-t_0; & M_c(t_0)<0 \\ 2(t^i_{1cd}(t_0)-t_0); & M_c(t_0)>0, c\neq d, \end{cases}, \tag{41}$$

$$c,d\in\{l,r\},\quad i\in\{1,2\},\quad t_0\in[0,T),\quad T=\frac{2\pi}{\omega}\quad ,$$

*by a Newton method, combined with a linear prolongation scheme. As the initial guesses for the Newton method we use the approximated values for $t^i_0$ and $\varepsilon^i_{cd}(t^i_0)$. See [20] for details on how this approximated values may be obtained.*

*When approximation to leading order in $H$ for $P(H)$ (as in equation (6)) is used to calculate $\varepsilon^i_{cd}(t^i_0)$, one gets that to leading order in $\beta$, $|\varepsilon^i_{cd}(t^i_0)| = |\varepsilon^i_{cc}(t^i_0)|$ for $\delta=0$. Actually, more accurate approximations (as using higher terms and Newton method) show that they are different.*